\documentclass[12pt]{amsart} 
\usepackage{amsfonts,latexsym} 
\usepackage{amsthm}
\usepackage{amssymb, amsmath}

\newtheorem{theorem}{Theorem}[section]

\newtheorem{lemma}[theorem]{Lemma}
\newtheorem{corollary}[theorem]{Corollary}

\theoremstyle{definition}
\newtheorem{defn}{Definition}[section]

\theoremstyle{remark}

\title{Representations with small $K$ types}
\author{Seung Won Lee}
\date{\today}

\usepackage{amssymb}
\usepackage{mathtools}
\usepackage{tikz}
\usetikzlibrary{chains}

\tikzset{node distance=2em, ch/.style={circle,draw,on chain,inner sep=2pt},chj/.style={ch,join},every path/.style={shorten >=4pt,shorten <=4pt},line width=1pt,baseline=-1ex}

\let\dlabel=\alabel

\newcommand{\dnode}[2][chj]{%
\node[#1,label={below:\dlabel{#2}}] {};
}

\newcommand{\dnodenj}[1]{%
\dnode[ch]{#1}
}

\newcommand{\dnodebr}[1]{%
\node[chj,label={below right:\dlabel{#1}}] {};
}

\newcommand{\dydots}{%
\node[chj,draw=none,inner sep=1pt] {\dots};
}

\begin{document}
\raggedbottom
\let\thefootnote\relax
\maketitle

\begin{abstract}
Let $\mathfrak{g}_{\mathbb{R}}$ be a split real, simple Lie algebra with complexification $\mathfrak{g}$.  Let $G_{\mathbb{C}}$ be the connected, simply connected Lie group with Lie algebra $\mathfrak{g}$, $G_{\mathbb{R}}$ the connected subgroup of $G_{\mathbb{C}}$ with Lie algebra $\mathfrak{g}_{\mathbb{R}}$, and $G$ a covering group of $G_{\mathbb{R}}$ with a maximal compact subgroup $K$.  A complete classification of ``small" $K$ types is derived via Clifford algebras, and an analog, $P^{\xi}$, of Kostant's $P^{\gamma}$ matrix is defined for a $K$ type ${\xi}$ of principal series admitting a small $K$ type.  For the connected, simply connected, split real forms of simple Lie types other than type $C_n$, a product formula for the determinant of $P^{\xi}$ over the rank one subgroups corresponding to the positive roots is proved.  We use these results to determine cyclicity of a small $K$ type of principal series in the closed Langlands chamber and irreducibility of the unitary principal series admitting a small $K$ type.
\end{abstract}

\newpage

\begin{section}{Introduction}

$\indent$For a linear, real reductive Lie group with a maximal compact subgroup $K$ and the complexification $\mathfrak{g}$ of its Lie algebra, the spherical principal series representations $I_{\nu}$ play a central role.  These are the principal series containing a trivial $K$ type.  Important work on the spherical principal series has been done by Kostant.  In [Kos], Kostant studies $I_{\nu}$ via the determinant of his $P^{\gamma}$ matrix defined for each $K$ type $\gamma$ that occurs in $I_{\nu}$.  Kostant determines this determinant by proving a product formula over the rank one subgroups corresponding to the reduced restricted roots and calculating the determinant in the rank one cases.  He proves several implications of the determinant formula; in particular, the cyclicity of the trivial $K$ type of $I_{\nu}$ in the closed Langlands chamber and the irreducibility of the unitary spherical principal series $I_{\nu} \ (Re \ \nu = 0)$.  This implies a complete classification of irreducible $(\mathfrak{g},K)$ modules admitting a trivial $K$ type as the unique irreducible quotients of $I_{\nu}$ in the closed Langlands chamber.

For a nonlinear, real reductive Lie group, the analogue of the spherical principal series representations are the genuine principal series representations $I_{\sigma,\nu}$ that admit a ``small'' $K$ type $(\tau,V_{\tau})$ defined by Wallach.  For the connected, simply connected, split real forms of simple Lie types other than type $C_n$, we study $I_{\sigma,\nu}$ via the determinant of a generalization of Kostant's $P^{\gamma}$ matrix, $P^{\xi}$, for a $K$ type $\xi$ that occurs in $I_{\sigma,\nu}$, where we determine this determinant by proving a similar product formula as Kostant's.  We omit type $C_n$ as the small $K$ types are $1$ dimensional, and representations that admit $1$ dimensional $K$ types are studied by Zhu [Zhu].

$\indent$ Let $G$ be a real reductive Lie group with a maximal compact subgroup $K$ and let $\mathfrak{g}_{\mathbb{R}} = Lie(G)$ and $\mathfrak{k}_{\mathbb{R}} = Lie(K)$.  We drop the subscripts to denote the complexifications of the Lie algebras.  Let $\mathfrak{g} = \mathfrak{k} \oplus \mathfrak{p}$ be a Cartan decomposition, let $\mathfrak{a}$ be a maximal abelian subspace of $\mathfrak{p}$ and let $^0\!M$ be the centralizer of $\mathfrak{a}$ in $K$.  Let $\mathcal{H}$ be the space of harmonics on $\mathfrak{p}$.  In 11.3.1 of [RRG II], Wallach defines a small $K$ type to be an irreducible representation $(\tau, V_{\tau})$ of $K$ such that $\sigma = V_{\tau}|_{{}^{0}\!M}$ is irreducible, and gives genuine examples for all simply connected, real forms of all simple Lie types.  Let $I_{P,\sigma,\nu}$ be the usual $K$ finite principal series for a minimal parabolic subgroup $P$, $\sigma$, and $\nu \in \mathfrak{a}^{*}$.  Let $Y^{\tau,\nu} = {U(\mathfrak{g}) \otimes_{U(\mathfrak{k})U(\mathfrak{g})^{K}} V_{\tau,\nu}}$ be the $(\mathfrak{g},K)$ module where $\mathfrak{g}$ acts by left translation, $K$ acts by $\rm{Ad} \otimes \tau$, and $U(\mathfrak{g})^{K}$ acts on $V_{\tau,\nu}$ as it does on $I_{P,\sigma,\nu}(\tau)$.  Wallach also proves that $\mathcal{H} \otimes V_{\tau} \cong Y^{\tau,\nu} \cong Ind_{^0\!M}^{K}(\sigma)$ as $K$ modules.  We also note that the spherical principal series is a special case as the trivial $K$ type is a small $K$ type.

Let now $\mathfrak{g}_{\mathbb{R}}$ be a split real, simple Lie algebra.  Let $G_{\mathbb{C}}$ be the connected, simply connected Lie group with Lie algebra $\mathfrak{g}$, $G_{\mathbb{R}}$ the connected subgroup of $G_{\mathbb{C}}$ with Lie algebra $\mathfrak{g}_{\mathbb{R}}$, and $G$ a covering group of $G_{\mathbb{R}}$ different from $G_{\mathbb{R}}$ with covering homomorphism $p$ and a maximal compact subgroup $K$.  From this point and on, a small $K$ type will always be a genuine representation of $K$ unless specified otherwise.  The first main result is the following classification of small $K$ types.  $Z = ker(p)$.  If $K$ is a product of two groups, denote by $p_j$ the projection onto the $j^{\rm{th}}$ factor for $j=1,2$.  Denote by $s$ the spin representation of $Spin(n)$ for $n$ odd, and either of the two half spin representations of $Spin(n)$ for $n$ even.

\vspace{.15 in}

\noindent \textbf{Theorem 1} {(Classification of Genuine Small $K$ types)}
\begin{center}
\begin{tabular}{ | l | p{4cm} | l | p{4cm} |}
    \hline
    $\rm{Type}$ & $K$ & $Z$ &  $\tau$  \\ \hline
    $A_n \ (n\geq 2)$ & $Spin(n+1)$ & $\mathbb{Z}/2\mathbb{Z}$  & $s$ \\ \hline
    $B_n \ (n\geq 3)$ & $Spin(n+1) \times Spin(n)$ & $\mathbb{Z}/2\mathbb{Z}$ & $s \circ p_1$ or $s \circ p_2$ for $n$ odd, $s \circ p_2$ for $n$ even \\ \hline
    $C_n$ & $SO(2) \times SU(n)$ & $\mathbb{Z}/m\mathbb{Z}$  & $\{e^{\frac{\pi i k}{m}} \ | \ (k,m)=1 \}$ \\ \hline
    $D_n \ (n\geq 3)$ & $Spin(n) \times Spin(n)$ & $\mathbb{Z}/2\mathbb{Z}$ & $s \circ p_1$ or $s \circ p_2$ \\ \hline
    $E_6 $ & $Sp(4)$ & $\mathbb{Z}/2\mathbb{Z}$  & $\mathbb{C}^{8}$ \\ \hline
    $E_7 $ & $SU(8)$ & $\mathbb{Z}/2\mathbb{Z}$  & $\mathbb{C}^{8}$ or $(\mathbb{C}^{8})^{*}$ \\ \hline
   $E_8 $ & $Spin(16)$ & $\mathbb{Z}/2\mathbb{Z}$  & $\mathbb{C}^{16}$ \\ \hline
    $F_4 $ & $Sp(3) \times SU(2)$ & $\mathbb{Z}/2\mathbb{Z}$  & $\mathbb{C}^{2}$ $\circ \ p_2$ \\ \hline
    $G_2 $ & $SU(2) \times SU(2)$ & $\mathbb{Z}/2\mathbb{Z}$  & $\mathbb{C}^{2} \circ p_1 \ or \ \mathbb{C}^{2} \circ p_2$ \\ \hline
\end{tabular}
\end{center}

\vspace{.15 in}

For a $K$ type ${\xi}$ that occurs in the principal series $I_{P,\sigma,\nu}$, a generalization of Kostant's $P^{\gamma}$ matrix, $P^{\xi}$, is defined (cf. Definition 3.1).  The definition of $P^{\xi}$ depends on the choice of a small $K$ type if $I_{P,\sigma,\nu}$ admits more than one (cf. section 3).  The second main result is a product formula for the determinant of $P^{\xi}$ over the rank one subgroups corresponding to the positive roots for $G$ other than type $C_n$.  Let $\phi$ be a positive root of $\mathfrak{g}$, and let $G_{\phi}$ be the corresponding rank one subgroup.  $G_{\phi}$ has its semisimple part the group generated by the metaplectic group $Mp_2(\mathbb{R})$ and $^0\!M$.  Let $K_{\phi}$ be the maximal compact subgroup of $G_{\phi}$ generated by a torus and $^0\!M$.  $V_{\xi} = \oplus_{j=1}^{n(\xi)} V_{\tau_j^{\phi}} \oplus U_{\phi}$ where $V_{\tau_j^{\phi}}$ is an irreducible $K_{\phi}$ submodule of $V_{\xi}$ such that $V_{\tau_j^{\phi}} \cong V_{\tau}$ as $^0\!M$ modules for all $j=1,...,n(\xi)$, with $n(\xi) = dim \ Hom_{{^0\!M}}(V_{\tau},V_{\xi})$.  Let $p_{\phi} =  p_{\tau_1}^{\phi}...p_{\tau_{n(\xi)}}^{\phi}$ where $p_{\tau_j}^{\phi}$ is the determinant of $P^{\tau_j^{\phi}}$ matrix of the rank one case of $G_{\phi}$ with $K_{\phi}$ type ${\tau_j^{\phi}}$.  Denote by $p_{(\phi)} = T_{\rho_{\phi}-\rho}(p_{\phi})$ where $T_{\rho_{\phi}-\rho}$ is translation by $\rho_{\phi}-\rho$.

\noindent \textbf{Theorem 2}  Let $p_{\xi}(\nu)$ denote the determinant of $P^{\xi}(\nu)$.  Then, there exists a nonzero scalar $c$ such that 
\[ p_{\xi}(\nu) = c\Pi_{\phi \in \Phi^{+}} p_{(\phi)}(\nu) \]

For the proof of the product formula we follow Kostant's strategy of proving that the rank one factors divide and are relatively prime, and proving that the degree is correct.  The divisibility argument is similar to Kostant's; however, our approach to the degree is different.  For this, we make a comparison of certain weight vectors of the torus corresponding to each positive root.

$\indent$ Theorem 2 has several implications analogous to the implications of Kostant's determinant formula.  In the simply laced cases and the doubly laced cases $B_n$ with the small $K$ type $s \circ p_2$, $F_4$ with the small $K$ type $\mathbb{C}^2 \circ p_2$, and $G_2$ with the small $K$ type $\mathbb{C}^2 \circ p_1$ ($p_1(K)$ is the short $SU(2)$), the implications are completely analogous to Kostant's.  First, cyclicity of the small $K$ type $I_{P,\sigma,\nu}(\tau)$ in the closed Langlands chamber (cf. Theorem 6.4) is proved as a consequence of the nonvanishing determinant of $P^{\xi}(\nu)$ in the closed Langlands chamber for all $\xi \in \hat{K}$ that occur in $I_{P,\sigma,\nu}$.  Second, irreducibility of the unitary principal series $I_{P,\sigma,\nu}$ $(Re \ \nu = 0)$ (cf. Theorem 6.5) is proved as a consequence of the cyclicity of $I_{P,\sigma,\nu}(\tau)$.  In the rest of the doubly laced cases, we note some differences in the implications.  For type $B_n$ with the small $K$ type $\tau = s \circ p_1$, $I_{P,\sigma,\nu}(\tau)$ is cyclic in the closed Langlands chamber for all other than $\nu$ s.t. $(\nu,\alpha) = 0$ for some $\alpha$ short (cf. Theorem 6.4).  This implies irreducibility of the unitary principal series $I_{P,\sigma,\nu}$ other than $\nu$ s.t. $(\nu,\alpha) = 0$ for some $\alpha$ short (cf. Theorem 6.5).  For type $G_2$ with the small $K$ types $\mathbb{C}^2 \circ p_j$ for $j=1, \ 2$, irreducibility of the unitary principal series still holds (cf. Theorem 6.5); however, cyclicity of the small $K$ type $\tau = \mathbb{C}^2 \circ p_2$ in the closed Langlands chamber does not hold.  In this case, $I_{P,\sigma,\nu}(\tau)$ is cyclic in the closed Langlands chamber if and only if $\frac{2(\nu,\alpha)}{(\alpha,\alpha)} \neq \frac{1}{2}$ for all $\alpha \in \Phi^{+}$ short (cf. Theorem 6.4).

\vspace{.1 in}

\noindent \textbf{Summary}
\begin{center}
\begin{tabular}{ | l | l | p{3.8 cm} | p{3.8 cm} |}
    \hline
    $\rm{Type}$ & $\tau$ & $\tau \rm{\ in \ the \ Closed}$ \text{Langlands Chamber} & $\rm{Unitary}$ \text{Principal Series} \\ \hline
    $\text{Simply Laced}$ & Any & Cyclic  & Irreducible \\ \hline
    $B_n \ (n\geq 3)$ & $s \circ p_1$ & Not always Cyclic & Sometimes reducible \\ \hline
    $B_n \ (n\geq 3)$ & $s \circ p_2$ & Cyclic & Irreducible \\ \hline
    $F_4$ & $\mathbb{C}^2 \circ p_2$ & Cyclic & Irreducible \\ \hline
    $G_2$ & $\mathbb{C}^2 \circ p_1$ & Cyclic & Irreducible \\ \hline
    $G_2$ & $\mathbb{C}^2 \circ p_2$ & Not always Cyclic & Irreducible \\ \hline
\end{tabular}
\end{center}

\vspace{.15 in}

In addition, a classification of irreducible $(\mathfrak{g},K)$ modules that admit a small $K$ type is derived, whose consequence for case $2$ below is nonexistence of discrete series representations that admit a small $K$ type $\tau$ (Corollary 6.7).

\vspace{.15 in}

\noindent \textbf{Theorem 3}  Let $G$ be a connected, simply connected, split real form of simple Lie type other than type $C_n$.  Let $V$ be an irreducible $(\mathfrak{g},K)$ module containing a small $K$ type $\tau$.
\begin{enumerate}
\item $V$ is $(\mathfrak{g},K)$ isomorphic with the unique irreducible quotient of $Y^{\tau,\nu}$ for some $\nu$ in the closed Langlands chamber.
\item If $\tau$ is any other than $s \circ p_1$ for type $B_n$ ($n$ necessarily odd) and $\mathbb{C}^2 \circ p_2$ for type $G_2$, $I_{P,\sigma,\nu} \cong Y^{\tau,\nu} = {U(\mathfrak{g}) \otimes_{U(\mathfrak{k})U(\mathfrak{g})^{K}} V_{\tau,\nu}}$ as $(\mathfrak{g},K)$ modules in the closed Langlands chamber.
\end{enumerate}

\vspace{.15 in}

The applications of the product formula are the following.  By the explicit computation of the determinant of the $P^{\xi}$ matrix in the rank one case (cf. 6.2) and the realization of the intertwining operator between the genuine principal series representations of $G$ as a ratio of the $P^{\xi}$ matrices (Theorem 6.2) similarly as in the spherical case (cf. [JW]), we derive a method to find the unknown shift factors of the determinant of the intertwining operators given in terms of ratios of the classical gamma functions derived by Vogan and Wallach [VW].  We derive explicit formulas for $\widetilde{SL(n,\mathbb{R})}$ (cf. 6.2).

Nonlinear covers of linear Lie groups and their representations have been studied in [RRG II] in general.  Our results are based on the exposition in Chapter 11 of [RRG II].  Properties of nonlinear, real groups are also studied in [A].  In split real cases, nonlinear groups and their representations are also studied in [AB], [ABPTV], [AT], [Hu], [PPS], and [Vog].  In [AB] and [PPS], genuine representations of the metaplectic groups are studied.  In [ABPTV], complementary series representations are studied where genuine representations of $^0\!M$ are discussed, and the formulas for the intertwining operators between the genuine principal series representations are derived.  In [AT], duality is studied for simply laced cases.  In [Hu], the unitary dual of the metalinear group is classified.  In [Vog], the unitary dual of the split $G_2$ is classified.  We note that the cyclicity assertions of the small $K$ type $\mathbb{C}^2 \circ p_2$ and the irreducibility of the genuine, unitary principal series are compatible with results in [Vog] (cf. Theorem 14.1 and Lemma 14.4 [Vog]).

We now discuss the layout of the article.  In section 2, we derive a complete classification of genuine, small $K$ types.  In section 3, we define for a $K$ type $\xi$ that occurs in the principal series representation admitting a small $K$ type $\tau$ a matrix $P^{\xi}$ (cf. Definition 3.1) that agrees with Kostant's $P^{\gamma}$ matrix [Kos] if $\tau$ is the trivial $K$ type.  In section 4, we derive structural results for $K$ modules that occur in $\mathcal{H}$ or $\mathcal{H} \otimes V_{\tau}$ where $\mathcal{H}$ is the space of harmonics on $\mathfrak{p}$.  In section 5, we prove a product formula for the determinant of $P^{\xi}$ over the rank one subgroups of $G$ corresponding to the positive roots.  In section 6, we derive implications and applications of the product formula for the determinant of $P^{\xi}$.

\vspace{.15 in}

\noindent \textbf{Acknowledgement}: The results stated in this article are based on the author's UC San Diego PhD Thesis (2012).  The author is grateful to his advisor, Professor Nolan Wallach, for the invaluable guidance and support throughout his graduate study.
 
\end{section}
\begin{section}{Classification of Small $K$ types}

Let $\mathfrak{g}$ be a simple Lie algebra over $\mathbb{C}$, and let $\mathfrak{g}_{\mathbb{R}}$ be the split real form of $\mathfrak{g}$.  Let $G_{\mathbb{C}}$ be the connected, simply connected Lie group with Lie algebra $\mathfrak{g}$, and let $G_{\mathbb{R}}$ be the connected subgroup of $G_{\mathbb{C}}$ with Lie algebra $\mathfrak{g}_{\mathbb{R}}$.  Let $G$ be a covering group of $G_{\mathbb{R}}$ different from $G_{\mathbb{R}}$ with covering homomorphism $p:G \rightarrow G_{\mathbb{R}}$.  Let $K$ be a maximal compact subgroup of $G$ defined as the subgroup of the fixed points of a Cartan involution $\theta$, and let $U$ be a compact form of $G_{\mathbb{C}}$ such that $G_{\mathbb{R}} \cap U = K_{\mathbb{R}} = p(K)$. 

In this section, we derive a complete classification of small $K$ types (see the definition below).  We first derive a complete classification of small $K$ types for $G$ of type $A_n$ via Clifford Algebras.  For other types, we introduce appropriate embeddings of $\widetilde{SL(n,\mathbb{R})}$ or the metalinear group $\widetilde{GL(n,\mathbb{R})}$ into $G$ to derive a complete classification of small $K$ types via the result for type $A_n$.

We also denote the differential of the above Cartan involution by $\theta$.  Let $\mathfrak{g}_{\mathbb{R}} = \mathfrak{k}_{\mathbb{R}} \oplus \mathfrak{p}_{\mathbb{R}}$ be the Cartan decomposition and $\mathfrak{a}_{\mathbb{R}}$ be a maximal abelian subspace of $\mathfrak{p}_{\mathbb{R}}$.  Let $^0\!M = {^0\!M_G} = Z_K(\mathfrak{a}_{\mathbb{R}})$ and ${^0\!M_{G_{\mathbb{R}}}} = Z_{K_{\mathbb{R}}}(\mathfrak{a}_{\mathbb{R}})$.

\begin{defn}{(Wallach)}
An irreducible representation $(\tau , V_{\tau})$ of $K$ is small if $\tau|_{^0\!M}$ is irreducible.
\end{defn}

Let $Z$ be the kernel of $p$.  $Z \subset K$.  Let $\chi$ be a genuine character of $Z$, i.e. $\chi$ is injective.  Theorem 1 in the introduction gives a complete classification of small $K$ types $(\tau, V_{\tau})$ such that $\tau|_{Z} = \chi Id$, i.e. genuine small $K$ types.  We assume throughout this article that $(\tau, V_{\tau})$ is a genuine small $K$ type unless stated otherwise.  We now prove Theorem 1 in the introduction.

\begin{proof}(of Theorem 1)

The $K$ types listed are genuine and small (cf. 11.A [RRG II]).  We prove the $K$ types give a complete list of small $K$ types.

As $G_{\mathbb{R}}$ is split and has a complexification $G_{\mathbb{C}}$, ${^0\!M_{G_{\mathbb{R}}} \cong (\mathbb{Z}/2\mathbb{Z})^{n}}$ where each $\mathbb{Z}/2\mathbb{Z}$ is from each and every simple root of $\mathfrak{g}_{\mathbb{R}}$ and $n$ is the rank of $\mathfrak{g}_{\mathbb{R}}$.

 For type $C_n$, as $Z \subset SO(2)$, $^0\!M = p^{-1}({^0\!M_{G_{\mathbb{R}}}})$ $= {p^{-1}((\mathbb{Z}/2\mathbb{Z})^{n})}$ is abelian.  Therefore, small $K$ types are necessarily 1 dimensional, and the $K$ types in the table exhaust the list. 

For the classification of small $K$ types for type $A_n \ (n\geq 2)$, we introduce Clifford Algebras.  Let $V = \mathbb{R}^{n+1}$, $( \ , \ )$ the standard inner product on $V$.  Let $T(V)$ be the tensor algebra on $V$, and let $I$ be the ideal of $T(V)$ generated by the elements $x\otimes x + (x,x)1$.  Set $Cliff(V) = T(V)/I$.  Let $\{ e_1,...,e_{n+1} \}$ be the standard basis of $V$.  $Cliff(V)$ is an algebra over $\mathbb{R}$ generated by $e_1,...,e_{n+1}$ with relations $e_j^2 = -1$ and $e_ie_j = - e_je_i$ if $i \neq j$.  Let $S = \{ v\in \mathbb{R}^{n+1} | (v,v)=1 \}$.  The maximal compact subgroup $Spin(n+1)$ of $\widetilde{SL(n+1,\mathbb{R})}$ is a subgroup of $Cliff(V)$ of products of even elements of $S$.  $Spin(n+1)$ is a two fold covering group of $SO(n+1)$, and $^0\!M_{\widetilde{SL(n+1,\mathbb{R})}}$ is generated by $\{e_ie_{i+1} | \ i=1,...,n \}$ (cf. 11.A.2.6 [RRG II]).  

$\mathbb{C}[{^0\!M_{\widetilde{SL(n+1,\mathbb{R})}}}]/ <{\eta+1}>$ is isomorphic to the Clifford Algebra on $\mathbb{C}^{n}$ where $\eta$ the nontrivial element of $Z = \mathbb{Z}/2\mathbb{Z}$.  $Cliff(\mathbb{C}^{n})$ is isomorphic to the simple matrix algebra $M_{2^{\frac{n}{2}}}(\mathbb{C})$ for $n$ even and a direct sum of two simple matrix algebras $M_{2^{\frac{n-1}{2}}}(\mathbb{C}) \oplus M_{2^{\frac{n-1}{2}}}(\mathbb{C})$ for $n$ odd (cf. Proposition 6.1.5 $\&$ 6.1.6 [GW]).  Hence a small $K$ type must be of dimension $2^{\frac{n}{2}}$ for $n$ even and $2^{\frac{n-1}{2}}$ for $n$ odd.  Weyl dimension formula implies the $K$ types in the table are the only ones with appropriate dimensions.

For types other than $A_n$, and $C_n$, we make use of the above analysis for type $A_n$.  We first realize certain embedding $i: S_{\mathbb{R}} \hookrightarrow G_{\mathbb{R}} = G/Z$ using Dynkin diagrams and extended Dynkin diagrams (cf. Chapter 6 [Bou]).

\begin{itemize}

\item $A_n$: $S_{\mathbb{R}} \cong GL(n,\mathbb{R})$.
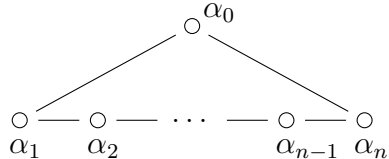
\begin{align*}
\begin{tikzpicture}[start chain]
\dnode{1}
\dnode{2}
\dydots
\dnode{n-1}
\dnode{n}
\begin{scope}[start chain=br going above]
\chainin(chain-3);
\node[ch,join=with chain-1,join=with chain-5,label={[inner sep=1pt]10:\( \alpha_0 \)}] {};
\end{scope}
\end{tikzpicture}
\end{align*}

Let $P$ be the parabolic subgroup with Levi factor $L$ where the simple roots of $Lie(L)$ are $\alpha_1,...,\alpha_{n-1}$.  The identity component of $L$ is isomorphic with $SL(n,\mathbb{R})$.  The embedded subgroup isomorphic to $GL(n,\mathbb{R})$ is generated by the $SL(n,\mathbb{R})$, ${\mathbb{Z}/2\mathbb{Z}} \subset {^0\!M_{G_{\mathbb{R}}}}$ from the node $\alpha_n$, and $\mathbb{R}_{>0}$ from the node $\alpha_0$.

\vspace{.1 in}

\item $B_n$ and $D_n$: $S_{\mathbb{R}} \cong SL(n,\mathbb{R})$.
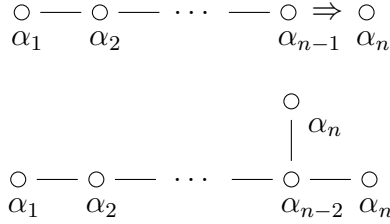
\begin{align*}
\begin{tikzpicture}[start chain]
\dnode{1}
\dnode{2}
\dydots
\dnode{n-1}
\dnodenj{n}
\path (chain-4) -- node[anchor=mid] {\(\Rightarrow\)} (chain-5);
\end{tikzpicture}
\end{align*}

\begin{align*}
\begin{tikzpicture}
\begin{scope}[start chain]
\dnode{1}
\dnode{2}
\node[chj,draw=none] {\dots};
\dnode{n-2}
\dnode{n-1}
\end{scope}
\begin{scope}[start chain=br going above]
\chainin(chain-4);
\dnodebr{n}
\end{scope}
\end{tikzpicture}
\end{align*}

Let $P$ be the parabolic subgroup with Levi factor $L$ where the simple roots of $Lie(L)$ are $\alpha_1,...,\alpha_{n-1}$.  The identity component of $L$ is isomorphic with $SL(n,\mathbb{R})$.

\vspace{.1 in}

\item $E_6$: $S_{\mathbb{R}} \cong GL(6,\mathbb{R})$.
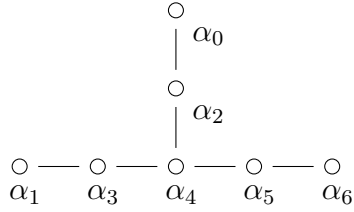
\begin{align*}
\begin{tikzpicture}
\begin{scope}[start chain]
\foreach \dyni in {1,3,4,5,6} {
\dnode{\dyni}
}
\end{scope}
\begin{scope}[start chain=br going above]
\chainin (chain-3);
\dnodebr{2}
\dnodebr{0}
\end{scope}
\end{tikzpicture}
\end{align*}

Let $P$ be the parabolic subgroup with Levi factor $L$ where the simple roots of $Lie(L)$ are $\alpha_1,\alpha_3,\alpha_4, \alpha_{5},\alpha_6$.  The identity component of $L$ is isomorphic with $SL(6,\mathbb{R})$.  The embedded subgroup isomorphic to $GL(6,\mathbb{R})$ is generated by the $SL(6,\mathbb{R})$, ${\mathbb{Z}/2\mathbb{Z}} \subset {^0\!M_{G_{\mathbb{R}}}}$ from the node $\alpha_2$, and $\mathbb{R}_{>0}$ from the node $\alpha_0$.

\vspace{.1 in}

\item $E_7$: $S_{\mathbb{R}} \cong SL(8,\mathbb{R})$.
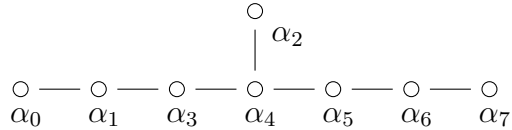
\begin{align*}
\begin{tikzpicture}
\begin{scope}[start chain]
\foreach \dyi in {0,1,3,4,5,6,7} {
\dnode{\dyi}
}
\end{scope}
\begin{scope}[start chain=br going above]
\chainin(chain-4);
\dnodebr{2}
\end{scope}
\end{tikzpicture}
\end{align*}

Let $P$ be the parabolic subgroup with Levi factor $L$ where the simple roots of $Lie(L)$ are $\alpha_0,\alpha_1,\alpha_3,\alpha_4, \alpha_{5},\alpha_6, \alpha_7$.  The identity component of $L$ is isomorphic with $SL(8,\mathbb{R})$.

\vspace{.1 in}

\item $E_8$: $S_{\mathbb{R}} \cong SL(9,\mathbb{R})$.
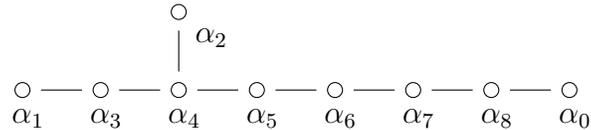
\begin{align*}
\begin{tikzpicture}
\begin{scope}[start chain]
\foreach \dyi in {1,3,4,5,6,7,8,0} {
\dnode{\dyi}
}
\end{scope}
\begin{scope}[start chain=br going above]
\chainin(chain-3);
\dnodebr{2}
\end{scope}
\end{tikzpicture}
\end{align*}

Let $P$ be the parabolic subgroup with Levi factor $L$ where the simple roots of $Lie(L)$ are $\alpha_1,\alpha_3,\alpha_4, \alpha_{5},\alpha_6, \alpha_7,\alpha_8, \alpha_0$.  The identity component of $L$ is isomorphic with $SL(9,\mathbb{R})$.

\vspace{.1 in}

\item $F_4$: $S_{\mathbb{R}} \cong Spin(5,4)_{0}$.
\begin{align*}
\begin{tikzpicture}[start chain]
\dnode{0}
\dnode{1}
\dnode{2}
\dnodenj{3}
\dnode{4}
\path (chain-3) -- node[anchor=mid]{\(\Rightarrow\)} (chain-4);
\end{tikzpicture}
\end{align*}

Let $P$ be the parabolic subgroup with Levi factor $L$ where the simple roots of $Lie(L)$ are $\alpha_0,\alpha_1,\alpha_2,\alpha_3$.  The identity component of $L$ is isomorphic with $Spin(5,4)_{0}$.

\vspace{.1 in}

\item $G_2$: $S_{\mathbb{R}} \cong SL(3,\mathbb{R})$.
\begin{align*}
\begin{tikzpicture}[start chain]
\dnode{1}
\dnodenj{2}
\dnode{0}
\path (chain-1) -- node{\(\Lleftarrow\)} (chain-2);
\end{tikzpicture}
\end{align*}

Let $P$ be the parabolic subgroup with Levi factor $L$ where the simple roots of $Lie(L)$ are $\alpha_2,\alpha_0$.  The identity component of $L$ is isomorphic with $SL(3,\mathbb{R})$.
\end{itemize}

The embedding $i: S_{\mathbb{R}} \hookrightarrow G_{\mathbb{R}}$ lifts to the embedding $\widetilde{i}: S \hookrightarrow G$ where $S$ is a subgroup of $G$ isomorphic to $p^{-1}(i(S_{\mathbb{R}}))$.  For types $B_n$ and $D_n$, $^0\!M \cong {^0\!M_{\widetilde{SL(n,\mathbb{R})}}} \times \mathbb{Z}/2\mathbb{Z}$ where $\widetilde{SL(n,\mathbb{R})} \cong S$ and the $\mathbb{Z}/2\mathbb{Z}$ can be either $(\pm 1,1)$ or $(1,\pm 1)$ where $\pm 1$ is the kernel of the covering homomorphism $p:Spin(m) \rightarrow SO(m)$ for $m = n$ or $n+1$ (cf. 11.A.2.8 [RRG II]).  For the rest of the types, $i({^0\!M_{S_{\mathbb{R}}}})={^0\!M_{G_{\mathbb{R}}}}$ as ${^0\!M_{G_{\mathbb{R}}} \cong (\mathbb{Z}/2\mathbb{Z})^{n}}$ where each $\mathbb{Z}/2\mathbb{Z}$ is from each and every simple root of $\mathfrak{g}_{\mathbb{R}}$ and $n$ is the rank of $\mathfrak{g}_{\mathbb{R}}$.  As ${^0\!M_{G}} = p^{-1}({^0\!M_{G_{\mathbb{R}}}}) = p^{-1}(i({^0\!M_{S_{\mathbb{R}}}})) = \widetilde{i}({^0\!M_{S}})$, we derive the following.

\begin{center}
    \begin{tabular}{ | l | p{2cm} | p{6cm} |}
    \hline
    $G$ & $S$ & $^0\!M_G$ isomorphic to  \\ \hline
    $A_n$ & $\widetilde{GL(n,\mathbb{R})}$ & $^0\!M_{\widetilde{GL(n,\mathbb{R})}}$ \\ \hline
    $E_6 $ & $\widetilde{GL(6,\mathbb{R})}$ & $^0\!M_{\widetilde{GL(6,\mathbb{R})}}$ \\ \hline
    $E_7 $ & $\widetilde{SL(8,\mathbb{R})}$ & $^0\!M_{\widetilde{SL(8,\mathbb{R})}}$ \\ \hline
    $E_8 $ & $\widetilde{SL(9,\mathbb{R})}$ & $^0\!M_{\widetilde{SL(9,\mathbb{R})}}$ \\ \hline
    $F_4 $ & $\widetilde{Spin(5,4)_{0}}$ & ${^0\!M_{\widetilde{Spin(5,4)_{0}}}} \cong {^0\!M_{\widetilde{SL(4,\mathbb{R})}} \times \mathbb{Z}/2\mathbb{Z}}$\\ \hline
    $G_2 $ & $\widetilde{SL(3,\mathbb{R})}$ & $^0\!M_{\widetilde{SL(3,\mathbb{R})}}$ \\ \hline
    \end{tabular}
\end{center}

We find necessary dimensions of small $K$ types using the discussion for type $A_n$.  Weyl dimension formula implies the $K$ types in the table exhaust the list of all small $K$ types.
\end{proof}

Consider now the metalinear group $\widetilde{GL(n,\mathbb{R})}$ and the group $\widetilde{Pin(n,n)}$ for $n\geq 3$ with maximal compact subgroups $Pin(n)$ and $Pin(n) \times Pin(n)$ respectively.  We give the following classification of small $K$ types for use in the proof of the product formula for $p_{\xi}$.  Let $s$ be the pin representation or either of the two pin representations of $Pin(n)$ depending on the parity of $n$.

\begin{theorem}
The small $K$ type for $\widetilde{GL(n,\mathbb{R})} \ (n\geq 3)$ is $s$ and the small $K$ types for $\widetilde{Pin(n,n)} \ (n\geq 3)$ are $s \circ p_1$ and $s \circ p_2$.

\begin{proof}
For $G$ the metalinear group $\widetilde{GL(n,\mathbb{R})}$, ${^0\!M_{\widetilde{GL(n,\mathbb{R})}}} \cong {^0\!M_{\widetilde{SL(n+1,\mathbb{R})}}}$ (cf. proof of Theorem 1).  Consider the embedding $i:O(n,n) \hookrightarrow SO(n+1,n+1)_{\circ}$ where the image of the maximal compact subgroup ${O(n) \times O(n)}$ of $O(n,n)$ under $i$ is contained in the maximal compact subgroup $SO(n+1) \times SO(n+1)$ of $SO(n+1,n+1)_{\circ}$ such that if $(g,h) \in O(n) \times O(n)$, 

\[ i((g,h)) = \begin{pmatrix}{}
det(g)^{-1} & 0 \\
0 & g \\
\end{pmatrix} , \begin{pmatrix}{}
det(h)^{-1} & 0 \\
0 & h \\
\end{pmatrix} \in SO(n+1) \times SO(n+1) \]  Let $p: \widetilde{Spin(n+1,n+1)_{\circ}} \rightarrow SO(n+1,n+1)_{\circ}$ be the covering homomorphism.  $p^{-1}(i(O(n,n)))$ is a Lie subgroup isomorphic to $\widetilde{Pin(n,n)}$, hence we have an embedding $\widetilde{i}:\widetilde{Pin(n,n)} \hookrightarrow \widetilde{Spin(n+1,n+1)_{\circ}}$.  \[ ^0\!M_{{SO(n+1,n+1)_{\circ}}} = \{ \begin{pmatrix}{}
det(g)^{-1} & 0 \\
0 & g \\
\end{pmatrix} , \begin{pmatrix}{}
det(g)^{-1} & 0 \\
0 & g \\
\end{pmatrix} | g \in O(n) \ diagonal  \} \] $= i(^0\!M_{{O(n,n)}})$.  As $\widetilde{i}(^0\!M_{\widetilde{Pin(n,n)}}) = p^{-1}(i(^0\!M_{{O(n,n)}}))$ and $^0\!M_{\widetilde{Spin(n+1,n+1)_{\circ}}} = p^{-1}(^0\!M_{{SO(n+1,n+1)_{\circ}}})$, $\widetilde{i}(^0\!M_{\widetilde{Pin(n,n)}}) = {^0\!M_{\widetilde{Spin(n+1,n+1)_{\circ}}}}$.  This in fact is true for $n\geq 2$.

The analysis for type $A_n$ gives the necessary dimensions of a small $K$ type $V_{\tau}$ and proves the $K$ types listed are the only small $K$ types.

\end{proof}
\end{theorem}

\end{section}

\begin{section}{$P^{\xi}$ Matrix}

Let $G$ be as in section 2.  In this section, we define a matrix, $P^{\xi}$, with entries in $U(\mathfrak{a})$ the universal enveloping algebra in $\mathfrak{a} = \mathfrak{a}_{\mathbb{R}} \otimes \mathbb{C}$ where $\xi$ is a $K$ type that occurs in the principal series representation of $G$ that admits a small $K$ type $\tau$.  If $\tau$ is the trivial $K$ type, the definition of $P^{\xi}$ matrix agrees with Kostant's $P^{\gamma}$ matrix [Kos].

Let $P = {^0\!M}AN$ be a parabolic subgroup of $G$ with a given Langlands decomposition.  Let $(\sigma,H_{\sigma})$ be an irreducible Hilbert representation of $^0\!M$ that is unitary when restricted to $K \cap {^0\!M}$.  Let $^{\infty}\!H^{P,\sigma,\nu}$ be the space of all smooth functions $f:G \longrightarrow H_{\sigma}$ such that $f(mang)=\sigma(m)a^{\nu}f(g)$ for $m \in {^0\!M}$, $a \in A$, $n \in N$, and $g \in G$ where $\nu \in Lie(A)_{\mathbb{C}}^{*}$.  Define for $f,g \in {^{\infty}\!H^{P,\sigma,\nu}}$ the inner product $\left \langle {f,g} \right \rangle = \int_{K} \left \langle f(k),g(k) \right \rangle_{\sigma} dk $.  Let $H^{P,\sigma,\nu}$ be the Hilbert space completion of $^{\infty}\!H^{P,\sigma,\nu}$.  The right regular action $\pi_{P,\sigma,\nu}(g)f(x) = f(xg)$ gives a Hilbert Representation $(\pi_{P,\sigma,\nu},H^{P,\sigma,\nu})$ of $G$ called the principal series representation.  If $X \in \mathfrak{g}_{\mathbb{R}}$, $X.f(g) = \frac{d}{dt}|_{t=0}f(g\cdot \rm{exp}(tX))$ defines a natural action of $\mathfrak{g}_{\mathbb{R}}$ on $H^{P,\sigma,\nu}$ induced from $\pi_{P,\sigma,\nu}$.  We also denote the action of $\mathfrak{g}_{\mathbb{R}}$ and hence the action of $\mathfrak{g}$ and its universal enveloping algebra $U(\mathfrak{g})$ by $\pi_{P,\sigma,\nu}$.  For $\gamma \in \hat{K}$, let $H^{P,\sigma,\nu}(\gamma)$ be the sum of all $K$ invariant, finite dimensional subspaces of $H^{P,\sigma,\nu}$ that are in the class of $\gamma$.  Let $I_{P,\sigma,\nu}$ be the algebraic direct sum $\oplus_{\gamma \in \hat{K}} H^{P,\sigma,\nu}(\gamma) \cap {^{\infty}\!H^{P,\sigma,\nu}}$, the underlying $(\mathfrak{g},K)$ module of the principal series representation $(\pi_{P,\sigma,\nu},H^{P,\sigma,\nu})$.

From this point on, let $P$ be a minimal parabolic subgroup of $G$ and let $\sigma = \tau|_{^0\!M}$ for a small $K$ type $\tau$.  We drop the subscript $\mathbb{R}$ to denote the complexifications of subspaces of $\mathfrak{g}_{\mathbb{R}}$ introduced in the last section.  Let $Y^{\tau,\nu}$ be the $(\mathfrak{g},K)$ module $U(\mathfrak{g})\otimes_{U(\mathfrak{g})^K U(\mathfrak{k})} V_{\tau,\nu}$ with the action of $\mathfrak{g}$ by left translation, the action of $K$ by $Ad \otimes \tau$, and the action of $U(\mathfrak{g})^K$ and $U(\mathfrak{k})$ on $V_{\tau,\nu}  = V_{\tau} \subset I_{P,\sigma,\nu}$ as a differential operator $\pi_{P,\sigma,\nu}$.  $Y_{\tau,\nu} \cong I_{P,\sigma,\nu}$ as $K$ modules (cf. 11.3.6 [RRG II]).  Let $\mathcal{H}$ be the space of harmonics on $\mathfrak{p}$, $\mathcal{J} = S(\mathfrak{p})^{K}$ the subspace of $K$ invariants of the space of symmetric polynomials on $\mathfrak{p}$.  Let $symm: S(\mathfrak{g}) \rightarrow U(\mathfrak{g})$ be the symmetrization map.  As $symm \otimes Id: S(\mathfrak{p}) \otimes U(\mathfrak{k}) \rightarrow U(\mathfrak{g})$ is a linear bijection and $symm(S(\mathfrak{p})) = symm(\mathcal{H}) \otimes symm(\mathcal{J})$, $symm(\mathcal{H}) \otimes V_{\tau} \cong I_{P,\sigma,\nu}$ as $K$ modules.

Let $(\xi, V_{\xi})$ be an irreducible representation of $K$ that occurs in $\mathcal{H} \otimes V_{\tau}$.  Let $n(\xi)$ be the multiplicity of $V_{\tau}|_{^0\!M}$ in $V_{\xi}$.  By Frobenius Reciprocity, the multiplicity of $\xi$ type in $\mathcal{H} \otimes V_{\tau}$ is $n(\xi)$.  Let $T^{\xi}_{1},...,T^{\xi}_{n(\xi)}$ be a basis of $Hom_{{}^{0}\!M}(V_{\tau},V_{\xi})$ and let $\epsilon^{\xi}_{1},...,\epsilon^{\xi}_{n(\xi)}$ be a basis of $Hom_{K}(V_{\xi},symm(\mathcal{H}) \otimes V_{\tau}) \cong Hom_K(V_{\xi},U(\mathfrak{g})\otimes_{U(\mathfrak{g})^K U(\mathfrak{k})} V_{\tau,\nu})$.  Let $\mu_{\tau,\nu}: U(\mathfrak{g})\otimes_{U(\mathfrak{g})^K U(\mathfrak{k})} V_{\tau,\nu} \longrightarrow I_{P,\sigma,\nu}$ be a $(\mathfrak{g},K)$ module homomorphism where the first factor acts on the second as a differential operator $\pi_{P,\sigma,\nu}$ (cf. 11.3.6 [RRG II]), and define $R_{\nu}:symm(\mathcal{H}) \otimes V_{\tau} \longrightarrow V_{\tau}$ by $R_{\nu}(Z):= \mu_{\tau,\nu}(Z)(e)$ where $e \in G$ is the identity element.  Every map defined in this paragraph intertwines ${}^{0}\!M$ action, hence $R_{\nu} \circ \epsilon_{i}^{\xi} \circ T_{j}^{\xi}$ also for all $i$ and $j$.

\newpage

\begin{defn} \
\begin{itemize}
\item Define by $P^{\xi}(\nu)$ the $n(\xi)$ by $n(\xi)$ matrix where $(P^{\xi}(\nu))_{i,j}$ is the polynomial in $\nu$ in which $R_{\nu} \circ \epsilon_{i}^{\xi} \circ T_{j}^{\xi}$ acts on $V_{\tau}$.  
\item Define by $P^{\xi}$ the $n(\xi)$ by $n(\xi)$ matrix $P^{\xi}(\nu)$ as an element of $U(\mathfrak{a}) = S(\mathfrak{a}) \cong \mathcal{O}(\mathfrak{a}^{*})$.  
\item  Define by $p_{\xi}$ and $p_{\xi}(\nu)$ the determinants of $P^{\xi}$ and $P^{\xi}(\nu)$ respectively.
\end{itemize}
\end{defn}	

\noindent $\mathbf{Remark \ 3.2}$  The definition of $P^{\xi}$ depends on the choice of a small $K$ type if $I_{P,\sigma,\nu}$ admits more than one small $K$ type.

\end{section}
\begin{section}{Structural Results}

Let $\mathfrak{g}_{\mathbb{R}}$ be any of the split real form of simple Lie type other than type $C_n$ or let $\mathfrak{g}_{\mathbb{R}} = \mathfrak{gl}(n,\mathbb{R})$ for $n \geq 3$ in which case $G$ is the metalinear group $\widetilde{GL(n,\mathbb{R})}$.  In this section, we derive structural results for certain $K$ modules.  In 4.1, we analyze the structure of $\mathcal{H}_{\alpha}$ and $\mathcal{H}_{\alpha} \otimes V_{\tau}$ for the rank one subgroup $G_{\alpha}$ of $G$ where $\alpha$ is a positive root of $\mathfrak{g}_{\mathbb{R}}$ and $\mathcal{H}_{\alpha}$ is $\mathcal{H}$ for $G_{\alpha}$.  In 4.2, we give an exposition of the irreducible $Pin(n)$ modules that can be found in 5.5.5 of [GW].  In 4.3, we prove Frobenius Reciprocity statements in our context for a $K$ module $V_{\xi}$ that occurs in $\mathcal{H} \otimes V_{\tau}$.

\begin{subsection}{}
Let $\mathfrak{g}_{\mathbb{R}} = \mathfrak{a}_{\mathbb{R}} \oplus \bigoplus_{\phi \in \Phi(\mathfrak{g}_{\mathbb{R}},\mathfrak{a}_{\mathbb{R}})} \mathfrak{g}_{\mathbb{R}}^{\phi}$ be the root space decomposition.  For $\alpha$ a positive root of $\mathfrak{g}_{\mathbb{R}}$, let ${\mathfrak{g}_{\mathbb{R}}}_{\alpha} = \mathfrak{a}_{\mathbb{R}} \oplus \mathfrak{g}_{\mathbb{R}}^{\alpha} \oplus \mathfrak{g}_{\mathbb{R}}^{-\alpha}$.  As $\mathfrak{g}_{\mathbb{R}}$ is split, $[{\mathfrak{g}_{\mathbb{R}}}_{\alpha},{\mathfrak{g}_{\mathbb{R}}}_{\alpha}] \cong \mathfrak{sl}(2,\mathbb{R})$.  Let $h_{\alpha} \in [{\mathfrak{g}_{\mathbb{R}}}_{\alpha},{\mathfrak{g}_{\mathbb{R}}}_{\alpha}]$ be such that $\alpha(h_{\alpha}) = 2$ and let $e_{\alpha} \in \mathfrak{g}_{\mathbb{R}}^{\alpha}$ be such that $[e_{\alpha},-\theta(e_\alpha)] = h_{\alpha}$.  $(h_{\alpha},e_{\alpha},-\theta(e_\alpha))$ is an $S$ triple.

\begin{defn}
Let $\alpha$ be a positive root of $\mathfrak{g}_{\mathbb{R}}$.  $t_{\alpha} := i(e_{\alpha}+\theta(e_{\alpha}))$.
\end{defn}

For a positive root $\alpha$ of $\mathfrak{g}_{\mathbb{R}}$, let $y_{\alpha} = e_{\alpha} - \theta(e_{\alpha})$, $Z_{\alpha} = h_{\alpha} + iy_{\alpha}$, and $\overline{Z_{\alpha}} = h_{\alpha} - iy_{\alpha}$.  Let $Z_{\alpha}^l$ and $\overline{Z_{\alpha}}^l$ be the $l^{th}$ tensor power of ${Z_{\alpha}}$ and $\overline{Z_{\alpha}}$ respectively for $l\in \mathbb{Z}_{\geq 0}$.  $t_{\alpha} \in \mathfrak{k}$, $[t_{\alpha},Z_{\alpha}^l] = -2lZ_{\alpha}^l$, and $[t_{\alpha},\overline{Z_{\alpha}}^l] = 2l\overline{Z_{\alpha}}^l$.  Let $\rm{exp}: \mathfrak{g}_{\mathbb{R}} \rightarrow G$ be the exponential map.  Let $G_{\alpha}$ be the rank one subgroup of $G$ generated by $^0\!M$ and the connected subgroup of $G$ with Lie algebra ${\mathfrak{g}_{\mathbb{R}}}_{\alpha}$.  Let $K_{\alpha}$ be the subgroup of $K$ generated by $\rm{exp}(i\mathbb{R}t_{\alpha})$ and ${}^{0}\!M$, the maximal compact subgroup of $G_{\alpha}$.  Let $\mathcal{H}_{\alpha}$ be the space of harmonics on $\mathfrak{p}_{\alpha} = \mathfrak{a} \oplus \mathbb{C}y_{\alpha}$ for the group $G_{\alpha}$.  $\mathcal{H}_{\alpha}$ as a space is $\oplus_{l\geq 0} \mathbb{C}Z_{\alpha}^l \bigoplus \oplus_{l>0} \mathbb{C}\overline{Z_{\alpha}}^l$.

\begin{lemma} Let $\alpha$ be a positive root of $\mathfrak{g}_{\mathbb{R}}$.  $\rm{Ad}({^0\!M})|_{t_{\alpha}} = \{ \pm 1 \}$.
\begin{proof}  As $Z \subset Z(G)$, it suffices to show $\rm{Ad}({^0\!M_{G_{\mathbb{R}}}})|_{t_{\alpha}} = \{ \pm 1 \}$.  First assume $\mathfrak{g}_{\mathbb{R}}$ is any other than $\mathfrak{gl}(n,\mathbb{R})$.  Since $G_{\mathbb{C}}$ is simply connected, the connected subgroup of $G_{\mathbb{C}}$ with Lie algebra $[\mathfrak{g}_{\beta},\mathfrak{g}_{\beta}]$ is isomorphic to $SL(2,\mathbb{C})$ for $\beta \in \Phi^{+}$ where $\mathfrak{g}_{\beta} = {\mathfrak{g}_{\mathbb{R}}}_{\beta} \otimes \mathbb{C}$.  Since $G_{\mathbb{R}}$ is split, the connected subgroup $G^{\beta}$ of $G_{\mathbb{R}}$ with Lie algebra $[{\mathfrak{g}_{\mathbb{R}}}_{\beta},{\mathfrak{g}_{\mathbb{R}}}_{\beta}]$ is isomorphic to $SL(2,\mathbb{R})$ for $\beta \in \Phi^{+}$.  Let $i_{\beta}:SL(2,\mathbb{R}) \rightarrow G^{\beta}$ be the isomorphism and let $\check{\beta}$ be the coroot to $\beta \in \Phi^{+}$.  $^0\!M_{G_{\mathbb{R}}}$ is generated by the elements $\rm{exp}(\pi i \check{\beta}) = i_{\beta}(-Id) \in G^{\beta}$ for $\beta \in \Phi^{+}$.  There exists a root $\beta$ such that $\alpha(\check{\beta})=1$.  For $\mathfrak{g}_{\mathbb{R}} = \mathfrak{gl}(n,\mathbb{R})$, $^0\!M_{G_{\mathbb{R}}} \cong {^0\!M_{{SL(n,\mathbb{R})}}} \times \mathbb{Z}/2\mathbb{Z}$.  Therefore, restriction to $\widetilde{SL(n,\mathbb{R})}$ proves the Lemma in this case.
\end{proof}
\end{lemma}

For $\mathfrak{g}_{\mathbb{R}}$, there are at most two root lengths, which we denote by short and long.  In the simply laced cases, we take the roots to be long.

\begin{lemma} Let $\alpha$ be a positive, long root of $\mathfrak{g}_{\mathbb{R}}$ and let $\tau$ be a small $K$ type.  The weights of $t_{\alpha}$ on $V_{\tau}$ are $\pm \frac{1}{2}$.
\begin{proof}

Consider first the cases $G = \widetilde{SL(n,\mathbb{R})}$ and $G = \widetilde{GL(n,\mathbb{R})}$ where $K=Spin(n)$ and $K=Pin(n)$ respectively.  For $K=Spin(n)$, $(\tau,V_{\tau})$ is the spin representation or either of the two half spin representations of $Spin(n)$ depending on the parity of $n$ by Theorem 1.  For $K=Pin(n)$, $(\tau,V_{\tau})$ is the pin representation or either of the two pin representations of $Pin(n)$ depending on the parity of $n$ by Theorem 2.1.  Hence the weights of $t_{\alpha}$ on $V_{\tau}$ are $\pm \frac{1}{2}$.

For $G$ other than $\widetilde{SL(n,\mathbb{R})}$ and $\widetilde{GL(n,\mathbb{R})}$, recall from the proof of Theorem 1 the embedding $\widetilde{i}:\widetilde{SL(n,\mathbb{R})} \ \hookrightarrow G$ or $\widetilde{i}:\widetilde{GL(n,\mathbb{R})} \ \hookrightarrow G$ for appropriate $n$ such that $\widetilde{i}(^0\!M_{\widetilde{SL(n,\mathbb{R})}}) \subseteq {^0\!M_{G}}$ or $\widetilde{i}(^0\!M_{\widetilde{GL(n,\mathbb{R})}}) \subseteq {^0\!M_{G}}$.  The containment is an equality for all other than the cases $B_n$, $D_n$, and $F_4$ where $^0\!M_G \cong {^0\!M_{\widetilde{SL(n,\mathbb{R})}}} \times \mathbb{Z}/2\mathbb{Z}$ for $B_n$ and $D_n$ and $^0\!M_G \cong {^0\!M_{\widetilde{SL(4,\mathbb{R})}}} \times \mathbb{Z}/2\mathbb{Z}$ for $F_4$.  The $\mathbb{Z}/2\mathbb{Z}$ can be chosen to act trivially on $V_{\tau}$ (cf. proof of Theorem 1).  In conclusion, $V_{\tau}$ restricted to either $\widetilde{i}(^0\!M_{\widetilde{SL(n,\mathbb{R})}})$ or $\widetilde{i}(^0\!M_{\widetilde{GL(n,\mathbb{R})}})$ is irreducible, hence $V_{\tau}$ is small as a $Spin(n)$ or a $Pin(n)$ module.  There is a simple, long root $\beta$ of $\mathfrak{g}_{\mathbb{R}}$ that restrict to a simple root of $Lie(\widetilde{i}(\widetilde{SL(n,\mathbb{R})}))$ or $Lie(\widetilde{i}(\widetilde{GL(n,\mathbb{R})}))$.  By the results for $\widetilde{SL(n,\mathbb{R})}$ and $\widetilde{GL(n,\mathbb{R})}$ above, the Lemma is proved for $\beta$.  As the positive, long roots are conjugates by elements of the Weyl group $W(A) = N_K(A)/Z_K(A)$ (cf. Proposition 6.11 [Bou]), the $t_{\alpha}$s are also (cf. Lemma 4.1).  Hence the Lemma is proved for any positive, long root $\alpha$ of $\mathfrak{g}_{\mathbb{R}}$.
\end{proof}
\end{lemma}

\begin{lemma} Let $\alpha$ be a positive, short root of $\mathfrak{g}_{\mathbb{R}}$.  $t_{\alpha}$ acts trivially on the small $K$ types $s \circ p_2$, $\mathbb{C}^2 \circ p_2$ for types $B_n$, $F_4$ respectively, and the weights of $t_{\alpha}$ are $\pm \frac{1}{2}$, $\pm \frac{3}{2}$, $\pm 1 $ on the small $K$ types $\mathbb{C}^2 \circ p_1$, $\mathbb{C}^2 \circ p_2$, $s \circ p_1$ for types $G_2$, $G_2$, and $B_n$ respectively.
\begin{proof}

For type $B_n$, $Lie(\mathfrak{k}) = \mathfrak{so}(n+1,\mathbb{C}) \oplus \mathfrak{so}(n,\mathbb{C})$.  Denote by $p_1$ the projection onto the first summand.  There is a positive, short root $\beta$ such that $p_1(t_{\beta}) = t_{\beta}$, with the only nonzero entries $(p_1(t_{\beta}))_{1,2} = 2i$ and $(p_1(t_{\beta}))_{2,1} = -2i$.  Thus the assertions are proved for $\beta$.

For type $F_4$, recall the embedding $\widetilde{i}: \widetilde{Spin(5,4)_0} \hookrightarrow G$ from the proof of Theorem 1.  The small $K$ type $\mathbb{C}^2 \circ p_2$ must restrict to the small $Spin(5) \times Spin(4)$ type $s \circ p_2$.  As the short, simple root $\beta$ for type $F_4$ restricts to a short, simple root for type $B_4$, the analysis for type $B_n$ above proves the assertion for $\beta$.

Consider now type $G_2$.  Let $\alpha$ be the short, simple root.  $\rm{exp}(i\pi t_{\alpha})^2 = -Id \times -Id \in SU(2) \times SU(2)$ (cf. 2.12 [Vog]).  Hence $\rm{exp}(i\pi t_{\alpha})$ must act on the two small $K$ types by $\pm \mathrm{exp}(\pi i (k+\frac{1}{2}))$ and $\pm \mathrm{exp}(\pi i (l+\frac{1}{2}))$ where $k,l \in \mathbb{Z}$.  As $^0\!M$ acts irreducibly on the small $K$ types and $\rm{Ad}({^0\!M})|_{t_{\alpha}} = \{ \pm 1 \}$ by Lemma 4.1, $t_{\alpha}$ weight on the two small $K$ types must be $\pm \frac{p}{2}$ and $\pm \frac{q}{2}$ respectively for some $p, \ q$ odd.

Consider the $7$ dimensional representation of the group $G$ of type $G_2$ with highest weight $[1,0]$ and weights $[2,-1], [-2,1], [\pm 1,0], [0,0], [1,-1]$, $[-1,1]$ in terms of the basis of fundamental weights.  $t_{\alpha}$ is the coroot to $\alpha$.  The restriction of the highest weight module $V_{[1,0]}$ to $K = SU(2) \times SU(2)$ ($p_1(K)$ and $p_2(K)$ are the short and long $SU(2)$s respectively) is $(\mathbb{C}^{2} \otimes \mathbb{C}^{2}) \oplus (\mathbb{C}^{3} \otimes \mathbb{C})$.  It is easily seen that the $t_{\alpha}$ weight of $0$ must correspond to a vector in $\mathbb{C}^3 \otimes \mathbb{C}$.  Let $X, \ Y$ be the dominant $t_{\alpha}$ weights on the small $K$ types $\mathbb{C}^{2} \circ p_1$ and $\mathbb{C}^{2} \circ p_2$ respectively.  If $X+Y =1$, then $X=Y=\frac{1}{2}$, hence $X-Y=0$, which does not occur in $\mathbb{C}^{2} \otimes \mathbb{C}^{2}$.  $X, \ Y$ thus must satisfy either $X+Y =2, \ X-Y =1$ or $X+Y =2 , \ X-Y = -1$.  In the former case, $X=\frac{3}{2}$ and $Y = \frac{1}{2}$, hence $t_{\alpha}$ weights on $\mathbb{C}^{3} \otimes \mathbb{C}$ must be $-3,0,3$.  But $\pm 3$ do not appear in the list of weights.  Therefore, the latter case must be true, $X = \frac{1}{2}$ and $Y=\frac{3}{2}$.

Similarly as in the proof of Lemma 4.2, the assertions are proved for any positive, short root $\alpha$.

\end{proof}
\end{lemma}

For a positive root $\alpha$ of $\mathfrak{g}_{\mathbb{R}}$, let $^0\!M_{\alpha}^{\pm}$ be the subsets of $^0\!M$ whose elements act on $t_{\alpha}$ by $\pm 1$ respectively.  If $t_{\alpha}$ acts nontrivially on $V_{\tau}$, let $V_{\tau,\alpha}^{\pm}$ be the subspaces of $V_{\tau}$ that consist of positive (resp. negative) $t_{\alpha}$ weight vectors. 

\begin{lemma} For a positive root $\alpha$ of $\mathfrak{g}_{\mathbb{R}}$, let $V_{{\gamma}_{\alpha}}$ be an irreducible nontrivial $K_{\alpha}$ module that occurs in $\mathcal{H}_{\alpha}$.  There exist exactly two $t_{\alpha}$ weights on $V_{{\gamma}_{\alpha}}$, and they are negatives of each other.
\begin{proof} Consider the $t_{\alpha}$ weight vector of weight $2l$, $\overline{Z_{\alpha}^{l}}$.  ${}^{0}\!M$ is a group of finite order that centralizes $\mathfrak{a}$ and ${}^{0}\!M_{\alpha}^{\pm}$ act by $\pm 1$ on $t_{\alpha}$ by Lemma 4.1.  Therefore, ${}^{0}\!M_{\alpha}^{+}$ fixes $\overline{Z_{\alpha}^{l}}$ and ${}^{0}\!M_{\alpha}^{-}$ moves $\overline{Z_{\alpha}^{l}}$ to $Z_{\alpha}^{l}$.  Hence the weights are $2l$ and $-2l$.
\end{proof}
\end{lemma}
		
\begin{lemma}For a positive root $\alpha$ of $\mathfrak{g}_{\mathbb{R}}$, let $V_{{\xi}_{\alpha}}$ be an irreducible $K_{\alpha}$ module that occurs in $\mathcal{H}_{\alpha} \otimes V_{\tau}$.  If $t_{\alpha}$ acts nontrivially on $V_{{\xi}_{\alpha}}$, there exist exactly two $t_{\alpha}$ weights on $V_{{\xi}_{\alpha}}$, and they are negatives of each other.
\begin{proof} Let $v \in V_{{\xi}_{\alpha}}$ be a $t_{\alpha}$ weight vector of weight $c=2l\pm r$ where $r= 0, \ \pm 1 ,\ \text{or} \ \pm \frac{p}{2}$ for $p$ odd by Lemma 4.2 and Lemma 4.3.  If $m \in {}^{0}\!M$, $t_{\alpha} \cdot m \cdot v$ = $m \cdot m^{-1} \cdot t_{\alpha} \cdot m \cdot m^{-1} \cdot m \cdot v$ = $\pm m \cdot t_{\alpha} \cdot v$ = $\pm c*m \cdot v$.  Since ${}^{0}\!M$ acts irreducibly on $V_{{\xi}_{\alpha}}$ the result follows.
\end{proof}
\end{lemma}

\begin{theorem} For a positive root $\alpha$ of $\mathfrak{g}_{\mathbb{R}}$, let $V_{{\gamma}_{\alpha}}$ be an irreducible $K_{\alpha}$ module that occurs in $\mathcal{H}_{\alpha}$.  $V_{{\gamma}_{\alpha}}$ as a space is the span of $\{\overline{Z_{\alpha}^{l}}, Z_{\alpha}^{l}\}$ for some l.  Moreover, ${}^{0}\!M$ invariant elements are $\mathbb{C}(\overline{Z_{\alpha}^{l}}+ Z_{\alpha}^{l})$.
\begin{proof} The statement follows from Lemma 4.4. 
\end{proof}
\end{theorem}

\begin{theorem} For a positive root $\alpha$ of $\mathfrak{g}_{\mathbb{R}}$, let $V_{{\xi}_{\alpha}}$ be an irreducible $K_{\alpha}$ module that occurs in $\mathcal{H}_{\alpha} \otimes V_{\tau}$.  If $t_{\alpha}$ acts nontrivially on $V_{\tau}$, $V_{{\xi}_{\alpha}}$ as a space is either $(\overline{Z_{\alpha}^{l}} \otimes V_{\tau,\alpha}^{+})  \oplus ({Z_{\alpha}^{l}} \otimes V_{\tau,\alpha}^{-})$ or $(\overline{Z_{\alpha}^{l}} \otimes V_{\tau,\alpha}^{-})  \oplus ({Z_{\alpha}^{l}} \otimes V_{\tau,\alpha}^{+})$ for some $l$.
\begin{proof} As $G_{\alpha}$ is a real reductive group of inner type and $V_{\tau}$ is a small $K_{\alpha}$ type, $\mathcal{H}_{\alpha} \otimes V_{\tau} \cong Ind_{^0\!M}^{K_{\alpha}}(\sigma)$ as $K_{\alpha}$ modules (cf. 11.3.6 [RRG II]).  Consider $\overline{Z_{\alpha}^{l}} \otimes v$ for some $l$ and $v \in V_{\tau,\alpha}^{+}$.  By Lemma 4.4 and Lemma 4.5, $K_{\alpha}$ module generated by this element is contained in $(\overline{Z_{\alpha}^{l}} \otimes V_{\tau,\alpha}^{+})  \oplus ({Z_{\alpha}^{l}} \otimes V_{\tau,\alpha}^{-})$.  Since an irreducible $K_{\alpha}$ module that occurs in $\mathcal{H}_{\alpha} \otimes V_{\tau}$ must be of dimension greater than or equal to $\mathrm{dim} \ V_{\tau}$, the inclusion is an equality.  We argue similarly if we start with a vector in $V_{\tau,\alpha}^{-}$.  Therefore, any irreducible $K_{\alpha}$ module that occurs in $\mathcal{H}_{\alpha} \otimes V_{\tau}$ has multiplicity one, and the statement follows.  \end{proof}
\end{theorem}

\begin{corollary} Let $\alpha$ be a positive root of $\mathfrak{g}_{\mathbb{R}}$.  There is a unique dominant $t_{\alpha}$ weight on an irreducible $K_{\alpha}$ module $V_{{\gamma}_{\alpha}}$ that occurs in $\mathcal{H}_{\alpha}$ and a unique dominant $t_{\alpha}$ weight on an irreducible $K_{\alpha}$ module $V_{{\xi}_{\alpha}}$ that occurs in $\mathcal{H}_{\alpha} \otimes V_{\tau}$.  
\begin{proof}
The two statements are direct consequences of Theorem 4.6 and Theorem 4.7.
\end{proof}
\end{corollary}

\end{subsection}

\begin{subsection}{} We give an exposition of irreducible $Pin(n)$ modules that can be found in 5.5.5 of [GW].  The exposition is given for the groups $O(n)$ and $SO(n)$ in [GW]; however, the exposition is also true for the groups $Pin(n)$ and $Spin(n)$.  If $n = 2k+1$ is odd, let $g_{0}= - Id \in O(2k+1)$.  If $n=2k$ is even, let $g_{0} \in O(2k)$ be the diagonal matrix whose entries are all $1$ except for last ${g_{0}}_{2k,2k}=-1$.  Let $p:Pin(n) \rightarrow O(n)$ be the covering homomorphism and let $\zeta$ be a choice of $p^{-1}(g_0)$.  Let $(\pi_{\lambda},V_{\lambda})$ be the irreducible representation of $Spin(n)$ with highest weight $\lambda$ and let $(\rho_{\lambda},V_{\lambda})$ be the induced representation $Ind_{Spin(2k)}^{Pin(2k)}(\pi_{\lambda})$.  Let $ \mathfrak{n}^{-} \oplus \mathfrak{h} \oplus \mathfrak{n}^{+}$ be the triangular decomposition of $Lie(Spin(n))_{\mathbb{C}}$.

\begin{theorem}{(Theorem 5.5.23 in [GW])} The irreducible regular representations of $Pin(2k+1)$ are of the form $(\pi_{\lambda}^{\epsilon},V_{\lambda}^{\epsilon})$, where $(\pi_{\lambda}^{\epsilon},V_{\lambda}^{\epsilon})$ restricted to $Spin(2k+1)$ is the highest weight representation $(\pi_{\lambda},V_{\lambda})$, and $\zeta$ acts on $V_{\lambda}^{\epsilon}$ by $\epsilon I$ where $\epsilon=\pm$. \end{theorem}

\begin{theorem}{(Theorem 5.5.24 in [GW])} Let $k\geq 2$.  The irreducible representation $(\sigma,W)$ of $Pin(2k)$ is one of the following two types.
\begin{itemize}
		\item Suppose $dim W^{\mathfrak{n}^{+}}=1$ and $\mathfrak{h}$ acts by the weight $\lambda$ on $W^{\mathfrak{n}^{+}}$. $(\sigma,W) \cong (\pi_{\lambda}^{\epsilon},V_{\lambda}^{\epsilon})$ where $(\pi_{\lambda}^{\epsilon},V_{\lambda}^{\epsilon})$ restricted to $Spin(2k)$ is the highest weight representation $(\pi_{\lambda},V_{\lambda})$, and $\zeta$ acts on $W^{\mathfrak{n}^{+}}$ by $\epsilon I$ where $\epsilon=\pm$.
		\item Suppose $dim W^{\mathfrak{n}^{+}}=2$.  Then $\mathfrak{h}$ has two distinct weights $\lambda$ and $\zeta \cdot \lambda$ on $W^{\mathfrak{n}^{+}}$, and $(\sigma,W) \cong (\rho_{\lambda},V_{\lambda})$.
\end{itemize}
\end{theorem}

\end{subsection}

\begin{subsection}{} For an irreducible representation $(\delta,V_{\delta})$ of a subgroup of $G$, let $d(\delta) = dim \ V_{\delta}$.  If $(\delta,V_{\delta})$ is an irreducible $K$ module that occurs in $\mathcal{H}$, let $l(\delta) = dim \ V_{\delta}^{^0\!M}$ and if $(\delta,V_{\delta})$ is an irreducible $K$ module that occurs in $\mathcal{H} \otimes V_{\tau}$, let $n(\delta) = dim \ Hom_K(V_{\delta},\mathcal{H} \otimes V_{\tau})$. 

\begin{lemma} Let ${\xi}$ be a $K$ type that occurs in $\mathcal{H} \otimes V_{\tau}$, and let ${\gamma_1},...,{\gamma_N}$ be the distinct $K$ types that occur in $\mathcal{H}$ such that $V_{\xi} \subset V_{\gamma_j} \otimes V_{\tau}$ for all $j=1,...,N$.  $n(\xi) = \Sigma_{j=1}^{N} l(\gamma_j)$.
\begin{proof}
For $G$ other than $\widetilde{GL(n,\mathbb{R})}$, $V_{\xi}$ has multiplicity one in $V_{\gamma_j} \otimes V_{\tau}$ as $V_{\tau}$ is multiplicity free (Proposition 3.2 [Ku]).  For $G = \widetilde{GL(n,\mathbb{R})}$, we restrict to $\widetilde{SL(n,\mathbb{R})}$ and use the branching of irreducible $Pin(n)$ modules to $Spin(n)$ (cf. 4.2) to deduce that $V_{\xi}$ has multiplicity one in $V_{\gamma_j} \otimes V_{\tau}$.  Each of $V_{\gamma_j}$ has multiplicity $l(\gamma_j)$ in $\mathcal{H}$ by Frobenius Reciprocity.  Hence $n(\xi) = \Sigma_{j=1}^{N} l(\gamma_j)$.
\end{proof}
\end{lemma}

\begin{lemma} Let $G$ be any of the connected, simply connected, split real forms of simple Lie types other than type $B_n \ (n \ \text{even})$, $C_n$, and $F_4$, or let $G$ be any of $\widetilde{GL(n,\mathbb{R})}$ and $\widetilde{Pin(n,n)}$ for $n \geq 3$.  If $V_{\xi}$ is an irreducible $K$ module that occurs in $\mathcal{H} \otimes V_{\tau}$, $V_{\xi}|_{^0\!M_{G}} = \bigoplus_{j=1}^{\frac{d(\xi)}{d(\tau)}} V_{\tau_j}$ where $V_{\tau_j} \cong V_{\tau}$ as $^0\!M_{G}$ modules.
\begin{proof}

For $G = \widetilde{SL(n,\mathbb{R})}$, let $\eta$ be the nontrivial element of $Z = \mathbb{Z}/2\mathbb{Z}$.  $\mathbb{C}[{^0\!M_{\widetilde{SL(n,\mathbb{R})}}}]/ <{\eta+1}>$ is isomorphic to the simple matrix algebra $M_{2^{\frac{n}{2}}}(\mathbb{C})$ for $n$ odd and a direct sum of two simple matrix algebras $M_{2^{\frac{n-1}{2}}}(\mathbb{C}) \oplus M_{2^{\frac{n-1}{2}}}(\mathbb{C})$ for $n$ even (cf. proof of Theorem 1).  Hence the Lemma is proved for $n$ odd.  If $n$ is even, $V_{{\tau}_{j}}$ is either of the two half spin representations restricted to ${}^{0}\!M_{\widetilde{SL(n,\mathbb{R})}}$.  Let $\omega$ be a choice of $p^{-1}(-Id)$ where $p:Spin(n)\rightarrow SO(n)$ is the covering homomorphism.  $\omega$  distinguishes the two representations as $^0\!M_{\widetilde{SL(n-1,\mathbb{R})}}$ does not.  $\omega$ acts trivially on $\mathcal{H}$ as it is central in $K$ and is an element of ${}^{0}\!M_{\widetilde{SL(n,\mathbb{R})}}$.  Therefore $\omega$ acts on the entire space $\mathcal{H} \otimes V_{\tau}$ as it does on $V_{\tau}$, hence the Lemma is proved for $n$ even.  For $G = \widetilde{GL(n,\mathbb{R})}$ the metalinear group, the Lemma is proved similarly as ${^0\!M_{\widetilde{GL(n,\mathbb{R})}}} \cong {^0\!M_{\widetilde{SL(n+1,\mathbb{R})}}}$.  

For $G = \widetilde{Spin(n,n)_{0}}$, ${^0\!M_{G}} \cong {^0\!M_{\widetilde{SL(n,\mathbb{R})}}} \times \mathbb{Z}/2\mathbb{Z}$ where the $\mathbb{Z}/2\mathbb{Z} \subset K$ and can be chosen to act trivially on $\mathcal{H}$ and $V_{\tau}$ (cf. proof of Theorem 1).  Therefore, $(\mathcal{H} \otimes V_{\tau})|_{^0\!M_G}$ decomposes in the same way as $(\mathcal{H} \otimes V_{\tau})|_{^0\!M_{\widetilde{SL(n,\mathbb{R})}}}$ where the factor $^0\!M_{\widetilde{SL(n,\mathbb{R})}}$ is embedded diagonally in $K = Spin(n) \times Spin(n)$.  If $n$ is odd, the Lemma is proved similarly as in the case of $\widetilde{SL(n,\mathbb{R})}$.  If $n$ is even, there are two possibilities for $V_{\tau}|_{^0\!M_{\widetilde{SL(n,\mathbb{R})}}}$.  Let $\omega$ be a choice of $p^{-1}(-Id \times -Id)$ where $p$ is the covering homomorphism $p:Spin(n) \times Spin(n) \rightarrow SO(n) \times SO(n)$.  $\omega$ distinguishes the two representations as ${^0\!M_{\widetilde{SL(n-1,\mathbb{R})}}}$ does not.  $\omega$ acts trivially on $\mathcal{H}$ as it is central in $K$ and is an element of $^0\!M_{\widetilde{Spin(n,n)_{0}}}$.  Therefore $\omega$ acts on the entire space $\mathcal{H} \otimes V_{\tau}$ as it does on $V_{\tau}$, hence the Lemma is proved for $\widetilde{Spin(n,n)_{0}}$.  For $G = \widetilde{Pin(n,n)}$, the Lemma is proved similarly as ${^0\!M_{\widetilde{Pin(n,n)}}} \cong {^0\!M_{\widetilde{Spin(n+1,n+1)_0}}}$ (cf. proof of Theorem 2.1).

For $G = \widetilde{Spin(n+1,n)_{0}}$, ${^0\!M_{G}} \cong {^0\!M_{\widetilde{SL(n,\mathbb{R})}}} \times \mathbb{Z}/2\mathbb{Z}$ where the $\mathbb{Z}/2\mathbb{Z} \subset K$ and can be chosen to act trivially on $\mathcal{H}$ and $V_{\tau}$ (cf. proof of Theorem 1).  Hence the decomposition of $V_{\xi}$ into ${^0\!M_{G}}$ modules only depend on the factor isomorphic to ${^0\!M_{\widetilde{SL(n,\mathbb{R})}}}$.  If $n$ is odd, the Lemma is proved similarly as above.

For $G$ of type $E_7$, $^0\!M$ is isomorphic to $^0\!M_{\widetilde{SL(8,\mathbb{R})}}$.  $Z(G) \cong \mathbb{Z}/4\mathbb{Z}$.  $Z(G) \subset Z(K)$ by Theorem 7.2.5 of [HAHS].  As $Z(K) = \{ zId \ | z^8 = 1 \}$, $Z(G)$ in $K$ consist of $\pm Id$ and $\pm i*Id$.  As $Z(G) \subset {^0\!M}$, $i*Id \in {^0\!M}$ acts trivially on $\mathcal{H}$.  $i*Id$ distinguishes the two ${^0\!M}$ irreducible representations $\mathbb{C}^{8}$ and $(\mathbb{C}^{8})^{*}$ as it acts by different signs on the two.  $i*Id$ acts on $\mathcal{H} \otimes V_{\tau}$ as it does on $V_{\tau}$.  Therefore, the Lemma is proved similarly as above.

For $G$ of type $E_6$, $E_8$, and $G_2$, $^0\!M$ is isomorphic to $^0\!M_{\widetilde{SL(n,\mathbb{R})}}$ for $n$ odd (cf. proof of Theorem 1).  Therefore the Lemma is proved similarly as the above case for type $\widetilde{SL(n,\mathbb{R})}$.
\end{proof} 
\end{lemma}

\begin{lemma} Let $G$ be any of the connected, simply connected, split real forms of simple Lie types other than type $C_n$.  Let $V_{\xi}$ be an irreducible $K$ module that occurs in $\mathcal{H} \otimes V_{\tau}$, $V_{\xi}|_{^0\!M_{G}} = \bigoplus_{j=1}^{n(\xi)} V_{\tau_j} \oplus U$ where $V_{\tau_j} \cong V_{\tau}$ as $^0\!M_{G}$ modules for $j=1,...,n(\xi)$.  Then, for any Weyl group element $s \in N_{K}(A)/Z_{K}(A)$, ${s \cdot V_{\tau_j}} \cong V_{\tau}$ as $^0\!M$ modules.
\begin{proof}
${s \cdot V_{\tau_j}}$ is an irreducible $^0\!M$ module as the Weyl group normalizes $^0\!M$ and the dimension is correct.  If $G$ is simply laced, of type $G_2$, or of type $B_{n}$ for $n$ odd, the Lemma is proved by Lemma 4.12.  Consider $G = \widetilde{Spin(n+1,n)_{0}}$ for $n$ even.  ${^0\!M_{G}} \cong {^0\!M_{\widetilde{SL(n,\mathbb{R})}}} \times \mathbb{Z}/2\mathbb{Z}$ where the $\mathbb{Z}/2\mathbb{Z} \subset K$ and can be chosen to act trivially on $\mathcal{H}$ and $V_{\tau}$.  Hence the decomposition of $V_{\xi}$ into ${^0\!M_{G}}$ modules only depend on the factor isomorphic to ${^0\!M_{\widetilde{SL(n,\mathbb{R})}}}$.  Let $\omega$ be the central element of $^0\!M_{\widetilde{SL(n,\mathbb{R})}}$ as in the proof of Lemma 4.12.  Note $s^{-1}\omega s$ is also central in $^0\!M_{\widetilde{SL(n,\mathbb{R})}}$.  $\omega$ distinguishes the two ${^0\!M_{\widetilde{SL(n,\mathbb{R})}}}$ representations.  The actions of $\omega$ and $s^{-1}\omega s$ on the small $K$ type $V_{\tau}$ are the same.  As $\omega s \cdot V_{\tau_j} = s s^{-1}\omega s \cdot V_{\tau_j}$, the Lemma is proved in this case.  For $G$ of type $F_4$, ${^0\!M_{G}} \cong {^0\!M_{\widetilde{SL(4,\mathbb{R})}}} \times \mathbb{Z}/2\mathbb{Z}$.  Choose $\omega$ similarly as above, and let $\upsilon$ be the nontrivial element of the factor $\mathbb{Z}/2\mathbb{Z}$.  The combination of the actions of $\omega$ and $\upsilon$ distinguish the four possible ${^0\!M_{G}}$ representations.  But, as $\omega$ and $\upsilon$ are central in ${^0\!M_{G}}$, the Lemma is proved similarly as in the case of $B_n$ above.
\end{proof} 
\end{lemma}

\end{subsection}

\end{section}
\begin{section}{Product Formula for $p_{\xi}$}
Let $\mathfrak{g}_{\mathbb{R}}$ be any of the split real forms of simple Lie types other than type $C_n$.  Let $G$ be the connected, simply connected Lie group with Lie algebra $\mathfrak{g}_{\mathbb{R}}$.  Recall the definition of $P^{\xi}$ for a $K$ type ${\xi}$ that occurs in $I_{P,\sigma,\nu}$.  The definition of $P^{\xi}$ depends on the choice of a small $K$ type if $I_{P,\sigma,\nu}$ admits more than one (cf. section 3).  In this section, we prove a product formula for $p_{\xi}$ over the rank one subgroups of $G$ corresponding to the positive roots.  The proof of the product formula consists of two parts, degree argument and divisibility argument.  In 5.1, we give a comparison of $t_{\alpha}$ weights on relevant vectors for $\alpha$ a positive root that is necessary for the degree argument.  In 5.2, we give the divisibility argument.  In 5.3, we prove a product formula for $p_{\xi}$.

Let $V_{\gamma}$ be an irreducible $K$ module that occurs in $\mathcal{H}$, and let $V_{\xi}$ be an irreducible $K$ module that occurs in $\mathcal{H} \otimes V_{\tau}$.  For a positive root $\alpha$ of $\mathfrak{g}_{\mathbb{R}}$, let $Span \ K_{\alpha} \cdot V_{\gamma}^{^0\!M} = \oplus_{j=1}^{l(\gamma)} W_{j}^{\alpha}$ be a decomposition into irreducible $K_{\alpha}$ submodules where $l(\gamma) = dim \ V_{\gamma}^{^0\!M}$ and $dim \ W_{j}^{\alpha} \cap V_{\gamma}^{^0\!M} = 1$ for $j=1,...,l(\gamma)$ (cf. Theorem 2.3.7 [Kos]).  Let $V_{\xi} = \oplus_{j=1}^{n(\xi)} V_{\tau_{j}^{\alpha}} \oplus U_{\alpha}$ be a decomposition into $K_{\alpha}$ submodules such that $V_{\tau_{j}^{\alpha}}$ is an irreducible $K_{\alpha}$ module with $V_{\tau_{j}}|_{^0\!M} \cong V_{\tau}$ for $j=1,...,n(\xi)$.  By Lemma 4.12, for $G$ other than type $B_n \ (n \ \mathrm{even})$ and $F_4$, $U_{\alpha} = 0$.

\begin{defn}
\
\begin{itemize}
\item Let $\delta_{\alpha,j}^{\gamma}$ be the unique dominant $t_{\alpha}$ weight on $W_{j}^{\alpha}$ for $j=1,...,l(\gamma)$ by Corollary 4.8.
\item Let $\delta_{\alpha,j}^{\xi}$ be the unique dominant $t_{\alpha}$ weight on $V_{\tau_j}^{\alpha}$ for $j=1,...,n(\xi)$ by Corollary 4.8.
\end{itemize}
\end{defn}

\begin{subsection}{}
Let $V_{\Gamma}$ be an irreducible $Pin(n)$ module that occurs in $\mathcal{H}$ the space of harmonics on $\mathfrak{p}$ for the metalinear group $\widetilde{GL(n,\mathbb{R})}$.  $\mathcal{H}$ is also that for $\widetilde{SL(n,\mathbb{R})}$.  Let $V_{\gamma}$ be an irreducible $Spin(n)$ module that occurs in $V_{\Gamma}|_{Spin(n)}$.  Recall the definition of $\zeta$ from 4.2.

\begin{theorem} Let $n=2k+1$.  $\zeta$ acts trivially on $V_{\Gamma}$, and there is no difference between $V_{\gamma}$ and $V_{\Gamma}$.
\begin{proof}
As the modules in the Theorem are submodules of $\mathcal{H}$, we may assume the modules are $SO(n)$ and $O(n)$ modules.  Recall the definition of $g_0$ from 4.2.  $g_0 \in Z(O(n))$ and $g_0 \in {}^{0}\!M_{{GL(n,\mathbb{R})}}$.  As $V_{\Gamma} \subset \mathcal{H}$ , $g_0$ must act trivially on $V_{\Gamma}$, hence there is no difference between $V_{\gamma}$ and $V_{\Gamma}$ by Theorem 4.9.
\end{proof}
\end{theorem}

\begin{theorem} Let $n=2k\geq 4$, and let $(m_1,...,m_k)$ be the highest weight of $V_{\gamma}$.
\begin{itemize}
\item If $m_k \neq 0$, $dim \ V_{\gamma}^{{}^{0}\!M_{\widetilde{SL(n,\mathbb{R})}}}$ = $dim \ V_{\Gamma}^{{}^{0}\!M_{\widetilde{GL(n,\mathbb{R})}}}$, $l(\gamma)=l(\Gamma)$ and $\delta_{\alpha,j}^{\gamma} = \delta_{\alpha,j}^{\Gamma}$ for all $j$ after reordering.
\item If $m_k = 0$, let $V_{\Gamma}$ = $V_{\gamma}^{\epsilon}$ where $\epsilon=\pm$ is the signature of $\zeta$ on the highest weight vector by Theorem 4.10.  $V_{\gamma}^{{}^{0}\!M_{\widetilde{SL(n,\mathbb{R})}}}$ = ${V_{\gamma}^{+}}^{{}^{0}\!M_{\widetilde{GL(n,\mathbb{R})}}} \oplus {V_{\gamma}^{-}}^{{}^{0}\!M_{\widetilde{GL(n,\mathbb{R})}}}$, hence $\delta_{\alpha,1}^{\gamma},...,\delta_{\alpha,l(\gamma)}^{\gamma}$ is a disjoint union of those of ${V_{\gamma}^{+}}^{{}^{0}\!M_{\widetilde{GL(n,\mathbb{R})}}}$ and ${V_{\gamma}^{-}}^{{}^{0}\!M_{\widetilde{GL(n,\mathbb{R})}}}$.
\end{itemize}
		
\begin{proof}
			As the modules in the Theorem are submodules of $\mathcal{H}$, we may assume the modules are $SO(n)$ and $O(n)$ modules.  Recall the definition of $g_0$ from 4.2.
			
			Let us assume first $m_k \neq 0$.  $g_{0}$ swaps the two $SO(2k)$ highest weight modules of highest weights $(m_1,...,m_k)$ and $(m_1,...,-m_k)$.  As $g_{0}$ commutes with ${}^{0}\!M_{SL(n,\mathbb{R})}$, $g_{0}$ gives a bijection of ${}^{0}\!M_{SL(n,\mathbb{R})}$-invariants in the $SO(2k)$ highest weight module of highest weight $(m_1,...,m_k)$ with ${}^{0}\!M_{SL(n,\mathbb{R})}$-invariants in the $SO(2k)$ highest weight module of highest weight $(m_1,...,-m_k)$.  Hence, $dim \ V_{\gamma}^{{}^{0}\!M_{{SL(n,\mathbb{R})}}}$ = $dim \ V_{\Gamma}^{{}^{0}\!M_{{GL(n,\mathbb{R})}}}$ as ${}^{0}\!M_{{GL(n,\mathbb{R})}}$ is generated by ${}^{0}\!M_{{SL(n,\mathbb{R})}}$ and $g_0$.  As $g_{0}$ leaves invariant $t_{\alpha}^2$, the statement of the $t_{\alpha}$ weights is also proved by Lemma 4.4.
			
			Assume $m_k = 0$.  Let $v(m_1,...,m_k)$ be the highest weight vector, and let $v_{1},...,v_{l(\gamma)}$ be a basis of $V_{\gamma}^{{}^{0}\!M_{{SL(n,\mathbb{R})}}}$ such that $g_{0}$ acts on $v_{j}$ by $\pm 1$ for all $j$, which is possible as $g_{0}^2=Id$ and $g_{0}$ commutes with ${}^{0}\!M_{SL(n,\mathbb{R})}$.  Denote by $v_{1}^{+},...,v_{l(\gamma)}^{+}$ and $v_{1}^{-},...,v_{l(\gamma)}^{-}$ the above basis thought of being in ${V_{\gamma}^{+}}$ and ${V_{\gamma}^{-}}$ respectively.  The only difference between ${V_{\gamma}^{+}}$ and ${V_{\gamma}^{-}}$ is the action of $g_{0}$.  Denote by $\epsilon$ the action of $g_{0}$ on $v_{1}^{+},...,v_{l(\gamma)}^{+}$.  $g_{0}$ acts by $-\epsilon$ on $v_{1}^{-},...,v_{l(\gamma)}^{-}$.  Indeed, $v_j = X_j.v(m_1,...,m_k)$ where $X_j\in U(\mathfrak{{n}}^{-})$.  As $g_{0}$ acts by different signatures on $v(m_1,...,m_k)$ for ${V_{\gamma}^{+}}$ and ${V_{\gamma}^{-}}$, the statement is now proved.
\end{proof}
\end{theorem}

Denote by $(\mathrm{T},V_{\mathrm{T}})$ the small $Pin(n)$ type and denote by $V_{\Xi}$ an irreducible $Pin(n)$ module that occurs in $\mathcal{H} \otimes V_{\mathrm{T}}$ for convenience.  Let $V_{\xi}$ be an irreducible $Spin(n)$ module that occurs in $V_{\Xi}|_{Spin(n)}$.  $V_{\xi} \subseteq  \mathcal{H} \otimes V_{\tau}$ where $V_{\tau}$ is a small $Spin(n)$ representation.  Let $V_{\xi}=\bigoplus_{j=1}^{n(\xi)} V_{\tau_{j}^{\alpha}}$ be a decomposition into irreducible $K_{\alpha}$ submodules such that $V_{\tau_{j}^{\alpha}} \cong V_{\tau}$ as $^0\!M_{\widetilde{SL(n,\mathbb{R})}}$ modules for all $j$ by Lemma 4.12.

\begin{theorem} If $n=2k\geq 4$, $V_{\Xi}|_{{}^{0}\!M_{\widetilde{GL(n,\mathbb{R})}}}$ = $\bigoplus_{j=1}^{n(\xi)} (V_{\tau_{j}^{\alpha}} \oplus \zeta \cdot V_{\tau_{j}^{\alpha}})$.  If $n=2k+1$, $V_{\Xi}|_{{}^{0}\!M_{\widetilde{GL(n,\mathbb{R})}}}$ = $\bigoplus_{j=1}^{n(\xi)} V_{\tau_{j}^{\alpha}}$.  In either case, $ \delta_{\alpha,j}^{\xi} = \delta_{\alpha,j}^{\Xi}$ for all $j$ after reordering.
\begin{proof}
As ${}^{0}\!M_{\widetilde{GL(n,\mathbb{R})}}$ is generated by ${}^{0}\!M_{\widetilde{SL(n,\mathbb{R})}}$ and $\zeta$, the Theorem is proved by Theorem 4.9, Theorem 4.10, and Lemma 4.12.
\end{proof}
\end{theorem}

Let ${\xi}$ be a $K$ type that occurs in $\mathcal{H} \otimes V_{\tau}$ and let ${\gamma_1},...,{\gamma_N}$ be the distinct $K$ types that occur in $\mathcal{H}$ such that $V_{\xi} \subset V_{\gamma_j} \otimes V_{\tau}$ for all $j=1,...,N$.  For $\alpha$ a positive root of $\mathfrak{g}_{\mathbb{R}}$, let $\delta_{\alpha,1}^{\xi},...,\delta_{\alpha,n(\xi)}^{\xi}$ be the $t_{\alpha}$ weights on $V_{\xi}$ and let $\delta_{\alpha,1},...,\delta_{\alpha,\Sigma_{j=1}^{N}l(\gamma_j)}$ be the $t_{\alpha}$ weights on $V_{\gamma_1},...,V_{\gamma_N}$ defined in Definition 5.1.  $n(\xi) = \Sigma_{j=1}^{N}l(\gamma_j)$ by Lemma 4.11.

\begin{theorem} If $\alpha$ is long, $\delta_{\alpha,1}^{\xi},...,\delta_{\alpha,n(\xi)}^{\xi}$ can be reordered so that $\delta_{\alpha,j} = \delta_{\alpha,j}^{\xi} \pm \frac{1}{2}$ for all $j=1,...,n(\xi)$.

\

\em{We will need the following Lemma for the proof of Theorem 5.4}.

\begin{lemma} The comparison of $t_{\alpha}$ weights in Theorem 5.4 for $\widetilde{SL(n,\mathbb{R})}$ implies that for $\widetilde{GL(n,\mathbb{R})}$.
\begin{proof}	
The space of Harmonics $\mathcal{H}$ on $\mathfrak{p}$ are the same for both $\widetilde{SL(n,\mathbb{R})}$ and $\widetilde{GL(n,\mathbb{R})}$.

Assume first $n$ is odd.  A small $Pin(n)$ representation $V_{\mathrm{T}}$ is the small $Spin(n)$ representation $V_{\tau}$ if we restrict from $Pin(n)$ to $Spin(n)$.  By Theorem 5.1, the restriction of $Pin(n)$ to $Spin(n)$ does not change the assumptions of the modules in Theorem 5.4.  Given an irreducible $Pin(n)$ module that occurs in $\mathcal{H}$, the $t_{\alpha}$ weights of interest do not change restricted to $Spin(n)$ by Theorem 5.1.  Given an irreducible $Pin(n)$ module that occurs in $\mathcal{H} \otimes V_{\mathrm{T}}$, the $t_{\alpha}$ weights of interest do not change restricted to $Spin(n)$ by Theorem 5.3.  Therefore, the comparison of the $t_{\alpha}$ weights for the group $Pin(n)$ reduces to the comparison of $t_{\alpha}$ weights for the group $Spin(n)$ by restriction of $Pin(n)$ to $Spin(n)$. 
		
Assume now $n$ is even.  The small $Pin(n)$ representation $V_{\mathrm{T}}$ is a direct sum of the two half spin representations $V_{\tau}$ and $\overline{V_{\tau}}$ if we restrict from $Pin(n)$ to $Spin(n)$.  Let ${\Gamma_1},...,{\Gamma_M}$ be the distinct $Pin(n)$ types that occur in $\mathcal{H}$ such that $V_{\Xi} \subseteq  V_{\Gamma_j} \otimes V_{\mathrm{T}}$.  Let us restrict $V_{\Gamma_j}$ to $Spin(n)$.  By Theorem 4.10, if $dim \ V_{\Gamma_j}^{\mathfrak{n}^{+}} = 2$, $V_{\Gamma_j}$ is a direct sum of two irreducible $Spin(n)$ modules with last entries of the highest weights nonzero and negatives of each other, and if $dim \ V_{\Gamma_j}^{\mathfrak{n}^{+}} = 1$, $V_{\Gamma_j}$ is irreducible as a $Spin(n)$ module.  For each $j$, let $V_{\gamma_j}$ be the choice of the irreducible $Spin(n)$ module that occurs in $V_{\Gamma_j}|_{Spin(n)}$ with last entry of the highest weight nonnegative, and reorder so that ${\gamma_1},...,{\gamma_N}$ are distinct $Spin(n)$ types.  $N\leq M$ as there may be $j$ such that $V_{\gamma_j}$ occurs twice with different $g_0$ signature on $V_{\gamma_j}^{\mathfrak{n}^+}$.  Let $V_{\xi}$ be the choice of the irreducible $Spin(n)$ module that occurs in $V_{\Xi}|_{Spin(n)}$ with last entry of the highest weight positive.  Without loss of generality, assume $V_{\xi} \subseteq \mathcal{H} \otimes V_{\tau}$.  ${\gamma_1},...,{\gamma_N}$ are distinct $Spin(n)$ types that occur in $\mathcal{H}$ such that $V_{\xi} \subseteq  V_{\gamma_j} \otimes V_{\tau}$.  Therefore, we can assume the statement of the $t_{\alpha}$ weights on these $Spin(n)$ modules.  But by Theorem 5.2 and Theorem 5.3, the comparison of $t_{\alpha}$ weights for the modules of the group $Pin(n)$ is that for $V_{\xi}, V_{\gamma_1},...,V_{\gamma_N}$ of $Spin(n)$.  Therefore, we have the result for $n$ even.
\end{proof}
\end{lemma}

\begin{proof}(of Theorem 5.4)

The $t_{\alpha}$ weights of interest are independent of the choice of $\alpha \in \Phi^{+}$.  Indeed, let $\psi$ and $\phi$ be two positive roots.  They are conjugates by an element of the Weyl group $N_K(A)/Z_K(A)$ hence $t_{\psi}$ and $t_\phi$ are also, by certain representative of the Weyl group element.  Let $s$ be such a representative.  $sK_{\psi}s^{-1} = K_{\phi}$.  First, $s \cdot V_{\gamma_j}^{^0\!M} = V_{\gamma_j}^{^0\!M}$
.  Now let $V_{\xi}=\oplus_{j=1}^{n(\xi)} V_{\tau_{j}^{\psi}} \oplus U_{\psi}$ be a decomposition into $K_{\psi}$ submodules such that each $V_{\tau_{j}^{\psi}}$ is an irreducible $K_{\psi}$ module isomorphic to $V_{\tau}|_{^0\!M}$ when restricted to ${^0\!M}$.  If $V_{\tau_{j}^{\phi}}=s \cdot V_{\tau_{j}^{\psi}}$ and $U_{\phi} = s \cdot U_{\psi}$, then $V_{\xi}=\oplus_{j=1}^{n(\xi)} V_{\tau_{j}^{\phi}} \oplus U_{\phi}$ is a decomposition into $K_{\phi}$ submodules such that each $V_{\tau_{j}^{\phi}}$ is an irreducible $K_{\phi}$ module isomorphic to $V_{\tau}|_{^0\!M}$ when restricted to ${^0\!M}$ by Lemma 4.13.
We first prove the Theorem for $\widetilde{SL(n,\mathbb{R})}$.

Let $n=3$.  Let $\xi=\frac{p}{2}$ be the highest weight of $V_{\xi}$ with $p$ odd.  If $p=1$, there exists only one $V_{\gamma} \subseteq \mathcal{H}$ with $V_{\xi} \subseteq V_{\gamma} \otimes V_{\tau}$ the trivial representation, and the claim is true.  If $p=3$, there exists only one $V_{\gamma} \subseteq \mathcal{H}$ with $V_{\xi} \subseteq V_{\gamma} \otimes V_{\tau}$ the representation with highest weight $2$.  In this case, the weights are $\frac{1}{2}$ and $\frac{3}{2}$ for $V_{\xi}$ and $0$ and $2$ for $V_{\gamma}$, hence the claim is also true.  Suppose $p>3$.  Then, there exist exactly two such representations, $V_{\gamma_{p_1}}$ and $V_{\gamma_{p_2}}$ with highest weights $\frac{p-1}{2}$ and $\frac{p+1}{2}$.  The weights of interest on $V_{\xi}$ are $\frac{1}{2}, \frac{3}{2},...,\frac{p}{2}$.  The weights of interest on $V_{\gamma_{p_1}}$ and $V_{\gamma_{p_2}}$ are $0,2,4,...,\frac{p-1}{2}$ and $2,4,...,\frac{p-1}{2}$ respectively if $\frac{p-1}{2}$ is even, and $2,4,...,\frac{p-3}{2}$ and $0,2,...,\frac{p+1}{2}$ respectively if $\frac{p+1}{2}$ is even as $^0\!M$ consist of elements $\pm Id, \pm \begin{pmatrix}{}
i & 0 \\
0 & -i \\
\end{pmatrix}, \pm \begin{pmatrix}{}
0 & 1 \\
-1 & 0 \\
\end{pmatrix}, \pm \begin{pmatrix}{}
0 & i \\
i & 0 \\
\end{pmatrix} $(cf. Lemma 8.11.8 [HAHS]).  In the first case where $\frac{p-1}{2}$ is even, consider $(\frac{1}{2}-\frac{1}{2}), (\frac{3}{2}+\frac{1}{2}), (\frac{5}{2}-\frac{1}{2}), (\frac{7}{2}+\frac{1}{2}),...,(\frac{p-2}{2}+\frac{1}{2}),(\frac{p}{2}-\frac{1}{2})$, which are $0,2,2,4,4,...,\frac{p-1}{2},\frac{p-1}{2}$.  This is exactly the union of the weights on $V_{\gamma_{p_1}}$ and $V_{\gamma_{p_2}}$.  In the second case where $\frac{p+1}{2}$ is even, consider $(\frac{1}{2}-\frac{1}{2}), (\frac{3}{2}+\frac{1}{2}), (\frac{5}{2}-\frac{1}{2}), (\frac{7}{2}+\frac{1}{2}),...,(\frac{p-2}{2}-\frac{1}{2}),(\frac{p}{2}+\frac{1}{2})$, which are $0,2,2,4,4,...,\frac{p-3}{2},\frac{p-3}{2},\frac{p+1}{2}$.  This is again exactly the union of the weights on $V_{\gamma_{p_1}}$ and $V_{\gamma_{p_2}}$.  Therefore the statement is proved for the case $n=3$.

We now proceed by induction.  Assume the Theorem for $\widetilde{SL(n,\mathbb{R})}$, and hence for $\widetilde{GL(n,\mathbb{R})}$ by Lemma 5.5.  We prove the Theorem for $\widetilde{SL(n+1,\mathbb{R})}$.  Recall from the proof of Theorem 1 the embedding $\widetilde{i}:\widetilde{GL(n,\mathbb{R})} \hookrightarrow \widetilde{SL(n+1,\mathbb{R})}$ with $\widetilde{i}(^0\!M_{\widetilde{GL(n,\mathbb{R})}}) = {^0\!M_{\widetilde{SL(n+1,\mathbb{R})}}}$.  We drop the notation $\widetilde{i}$ for simplicity.  We can restate the condition $V_{\xi} \subseteq V_{\gamma_j} \otimes V_{\tau}$ with $V_{\gamma_j} \subseteq V_{\xi} \otimes V_{\tau}^{*}$ where $V_{\tau}^{*}$ is the contragredient representation.  The statement of the Theorem is true with the restated condition for $\widetilde{GL(n,\mathbb{R})}$ by the induction hypothesis.

Let ${\gamma_1},...,{\gamma_N}$ be the distinct $Spin(n+1)$ types that occur in $\mathcal{H}$ such that $V_{\gamma_{j}} \subseteq V_{\xi} \otimes V_{\tau}^{*}$, and let $\bigoplus_{j=1}^{N} Span \ Pin(n) \cdot V_{\gamma_j}^{^0\!M} = \bigoplus_k W_{k}$ where each $W_{k}$ is an irreducible $Pin(n)$ module.  $V_{\xi}|_{Pin(n)} = \bigoplus_{j} V_{\xi_j}$ where each $V_{\xi_j}$ is an irreducible $Pin(n)$ submodule that occurs in $\mathcal{H} \otimes V_{\tau}$ by Lemma 4.12, with $\mathcal{H}$ that of $\widetilde{GL(n,\mathbb{R})}$.  $\bigoplus_k W_{k} \subseteq \bigoplus_{j} V_{\xi_j} \otimes V_{\tau}^{*}$ where each $W_k$ occurs in $\mathcal{H}$ of $\widetilde{GL(n,\mathbb{R})}$ as ${}^{0}\!M_{\widetilde{GL(n,\mathbb{R})}}$ = ${}^{0}\!M_{\widetilde{SL(n+1,\mathbb{R})}}$.  $V_{\xi_j} \otimes V_{\tau}^{*}$ decomposes into distinct $Pin(n)$ types as $V_{\tau}^{*}$ is multiplicity free (cf. proof of Lemma 4.11).  Therefore, if $W_k \cong W_l$ with $k \neq l$, $W_k$ and $W_l$ cannot be contained in a single $V_{\xi_j} \otimes V_{\tau}^{*}$.  This is important as the statement of the Theorem for $\widetilde{GL(n,\mathbb{R})}$ also assumes distinct ${\Gamma}$ types.  As the statement of the Theorem is true for $\widetilde{GL(n,\mathbb{R})}$ with the restated condition and the $t_{\alpha}$ weights of interest are the same after branching down to $Pin(n)$, the Theorem is proved for all positive roots $\alpha$ of $Lie(\widetilde{SL(n+1,\mathbb{R})})$ that restrict to a positive root of $Lie(\widetilde{GL(n,\mathbb{R})}) \subseteq Lie(\widetilde{SL(n+1,\mathbb{R})})$, hence for all positive roots of $Lie(\widetilde{SL(n+1,\mathbb{R})})$.

Let $G$ be other than $\widetilde{SL(n,\mathbb{R})}$.  Recall from the proof of Theorem 1 the embedding $\widetilde{i}:S \hookrightarrow G$ with $\widetilde{i}(^0\!M_{S}) \subseteq {^0\!M_G}$ where $S \cong \widetilde{SL(n,\mathbb{R})}$ or $\widetilde{GL(n,\mathbb{R})}$ for appropriate $n$.  Denote by $K_S$ the maximal compact subgroup of $S$.  We drop the notation $\widetilde{i}$ for simplicity.  $^0\!M_{S} = {^0\!M_G}$ for all but the cases of $B_n$, $D_n$, and $F_4$ where $^0\!M_{G} \cong {^0\!M_{\widetilde{SL(n,\mathbb{R})}}} \times \mathbb{Z}/2\mathbb{Z}$ for types $B_n$ and $D_n$, and $^0\!M_{G} \cong {^0\!M_{\widetilde{SL(4,\mathbb{R})}}} \times \mathbb{Z}/2\mathbb{Z}$ for type $F_4$.  The $\mathbb{Z}/2\mathbb{Z}$ can be chosen to act trivially on the small $K$ type.  Therefore, the small $K$ representation $V_{\tau}$ is also a small $K_S$ representation.  For $B_n$ and $D_n$, the $\mathbb{Z}/2\mathbb{Z} \subset Z(K)$ hence acts trivially on $\mathcal{H}$.

We can restate the condition $V_{\xi} \subseteq V_{\gamma_j} \otimes V_{\tau}$ as $V_{\gamma_j} \subseteq V_{\xi} \otimes V_{\tau}^{*}$.  The statement of the Theorem is true with the restated condition for $S$.  Let ${\gamma_1},...,{\gamma_N}$ be the distinct $K$ types that occur in $\mathcal{H}$ such that $V_{\gamma_{j}} \subseteq V_{\xi} \otimes V_{\tau}^{*}$, and let $\bigoplus_{j=1}^{N} Span \ K_{S} \cdot V_{\gamma_j}^{^0\!M_G} = \bigoplus_k W_{k}$ where each $W_{k}$ is an irreducible $K_{S}$ submodule that occurs in $\mathcal{H}$ of $S$ as $^0\!M_{G} = {^0\!M_S}$ or $^0\!M_{G} = {^0\!M_S} \times \mathbb{Z}/2\mathbb{Z}$.  $V_{\xi}|_{K_{S}} = \bigoplus_{j} V_{\xi_j} \oplus U_{\xi}$ where each $V_{\xi_j}$ is an irreducible $K_S$ module that occurs in $\mathcal{H} \otimes V_{\tau}$, with $\mathcal{H}$ that of $S$.  $\bigoplus_k W_{k} \subseteq \bigoplus_{j} V_{\xi_j} \otimes V_{\tau}^{*}$.  Indeed, if $W_{k} \subseteq U_{\xi} \otimes V_{\tau}^{*}$, then $W_{k} \otimes V_{\tau} \supseteq U_{\xi}$, a contradiction by Lemma 4.12.  For type $F_4$, $^0\!M_{G} \cong {^0\!M_{\widetilde{SL(4,\mathbb{R})}}} \times \mathbb{Z}/2\mathbb{Z}$ where the $\mathbb{Z}/2\mathbb{Z}$ centralizes $K_S$ (cf. proof of Theorem 1).  As the $\mathbb{Z}/2\mathbb{Z}$ acts trivially on $V_{\tau}$ and $W_k$, if $W_{k} \subseteq V_{\xi_j} \otimes V_{\tau}^{*}$, it must also act trivially on $V_{\xi_j}$.  Above observations are important as $\{ \delta_{\alpha, j}^{\xi} \}$ are from $^0\!M_G$ subrepresentations of $V_{\xi}$ isomorphic to $V_{\tau}|_{^0\!M_G}$.  $V_{\xi_j} \otimes V_{\tau}^{*}$ decomposes into distinct $K_S$ types as $V_{\tau}^{*}$ is multiplicity free (cf. proof of Lemma 4.11).  Therefore, if $W_k \cong W_l$ with $k \neq l$, $W_k$ and $W_l$ cannot be contained in a single $V_{\xi_j} \otimes V_{\tau}^{*}$.  This is important as the statement of the Theorem for $S$ also assumes distinct ${\gamma}$ types.  As the Theorem is true for $S$ with the restated condition and as the $t_{\alpha}$ weights of interest are the same on the appropriate modules after restricting to $K_S$, the Theorem is proved for all positive, long roots $\alpha$ of $Lie(G)$ that restrict to a positive, long root of $Lie(S) \subseteq Lie(G)$, hence for all positive, long roots of $Lie(G)$.
\end{proof}
\end{theorem}

Let $\alpha$ be a positive, short root of $\mathfrak{g}_{\mathbb{R}}$, and let $r$ be the unique dominant $t_{\alpha}$ weight on $V_{\tau}$ (cf. Lemma 4.3).  $n(\xi) = \Sigma_{j=1}^{N}l(\gamma_j)$ by Lemma 4.11.

\begin{theorem} Let $\alpha$ be a positive, short root of $\mathfrak{g}_{\mathbb{R}}$.  $\delta_{\alpha,j} = \delta_{\alpha,j}^{\xi} \pm r$ for all $j=1,...,n(\xi)$ after reordering.
\begin{proof}
Similarly as before, the $t_{\alpha}$ weights of interest are independent of the choice of the short root.

First consider $G$ of type $G_2$.  The maximal compact subgroup $K$ is $SU(2) \times SU(2)$.  $^0\!M_G \cong {^0\!M_{\widetilde{SL(3,\mathbb{R})}}}$ is embedded diagonally in $SU(2) \times SU(2)$ as $SO(3)$ must be embedded diagonally in $K_{\mathbb{R}} = SU(2) \times_{\mathbb{Z}/2\mathbb{Z}} SU(2)$ where $\mathbb{Z}/2\mathbb{Z} = \{(1,1),(-1,-1) \}$.  Note $p_1(K)$ is the short $SU(2)$ and $p_2(K)$ is the long $SU(2)$.  We also note that from the analysis of the case of $\widetilde{SL(3,\mathbb{R})}$ with $K=SU(2)$, $t_{\alpha}$ weight vectors of positive, even weight always arise from a $^0\!M$ invariant vector; however, a $t_\alpha$ weight vector of weight $0$ does not always arise from a $^0\!M$ invariant vector (cf. proof of Theorem 5.4).  Assume first $V_{\tau} = \mathbb{C}^2 \circ p_1$.

Let $V_r$ and $W_s$ be the irreducible representations of each of the two $SU(2)$s with highest weights $r$ and $s$ respectively.  The weights of $V_r$ are $-r,-r+1,...,r-1,r$ and similarly for $W_s$.  Denote these weight vectors by $v_{-r},v_{-r+1},...,v_{r-1},v_{r}$ and similarly for $W_s$.  Let $V_{\xi} \subseteq \mathcal{H} \otimes V_{\tau}$ with $V_{\xi} \cong V_r \otimes W_s$.  There is at most two distinct $V_{\gamma} \subseteq \mathcal{H}$ such that $V_{\xi} \subseteq  V_{\gamma} \otimes V_{\tau}$, $V_{\gamma_1} = V_{r_1} \otimes W_{s}$ and $V_{\gamma_2} = V_{r_2} \otimes W_{s}$ with $r_1 = r+\frac{1}{2}$ and $r_2 = r-\frac{1}{2}$.  Consider $v_i \otimes w_j \in V_{r_1} \otimes W_{s}$ and $v_i \otimes w_j \in V_{r_2} \otimes W_{s}$.  $t_{\alpha}$ weight on $v_i \otimes w_j$ is $i+3j$ by Lemma 4.3.  In order for them to be dominant $t_{\alpha}$ weight vectors that arise from $^0\!M$ invariant vectors, $i+3j$ is necessarily even and $i+3j \geq 0$ (cf. Theorem 4.6).  

First, if $r_1+3j$ is even and $r_1+3j>0$, then cover $v_{r_1} \otimes w_{j}$ with $v_{r} \otimes w_j$.  If $-r_1+3j>0$ and even also, cover $v_{-r_1} \otimes w_{j}$ with $v_{-r} \otimes w_{j}$.  Now assume $i+3j$ is even and $i+3j > 0$ with $i \neq r_1$ and $i \neq -r_1$.  We cover the two $v_i \otimes w_j \in V_{r_1} \otimes W_{s}$ and $v_i \otimes w_j \in V_{r_2} \otimes W_{s}$ with $v_{i\pm \frac{1}{2}} \otimes w_{j} \in V_r \otimes W_s$.  

Now assume $i+3j=0$ with $i,j \neq 0 $, and assume $i>0$.  In this case, $v_{i} \otimes w_{j} + v_{-i} \otimes w_{-j}$ is $^0\!M$ invariant, by scaling $v_{i} \otimes w_{j}$ appropriately if necessary (cf. Theorem 4.6).  If $i=r_1$, then cover the above with $v_{-i+\frac{1}{2}} \otimes w_{-j}$.  If $i\neq r_1$ then we cover the two $v_{i} \otimes w_{j} + v_{-i} \otimes w_{-j} \in V_{r_1} \otimes W_{s}$ and $v_{i} \otimes w_{j} + v_{-i} \otimes w_{-j} \in V_{r_2} \otimes W_{s}$ with $v_{i+\frac{1}{2}} \otimes w_{j}$ and $v_{-i+\frac{1}{2}} \otimes w_{-j}$.  The only ambiguity is when $i=j=0$.  But, as $l(V_{r_1}\otimes W_{s}) + l(V_{r_2}\otimes W_{s}) = n(V_{r}\otimes W_{s})$ by Lemma 4.11, and from the analysis of the case of $\widetilde{SL(3,\mathbb{R})}$ with $K=SU(2)$ (cf. proof of Theorem 5.4), one of the two $v_{0} \otimes w_{0} \in V_{r_1} \otimes W_{s}$ and $v_{0} \otimes w_{0} \in V_{r_2} \otimes W_{s}$ won't be $^0\!M$ invariant.  The case of $V_{\tau} = \mathbb{C}^2 \circ p_2$ is proved similarly.

Let now $G$ be of type $B_n$ with a maximal compact subgroup $K = Spin(n+1) \times Spin(n)$.  There is a short, positive root $\alpha$ such that $p_1(t_{\alpha}) = t_{\alpha}$, with the only nonzero entries $(p_1(t_{\alpha}))_{1,2} = 2i$ and $(p_1(t_{\alpha}))_{2,1} = -2i$ where $p_1$ denotes the projection of $\mathfrak{k}$ onto its first summand.  Let $S = \widetilde{Pin(n-1,n-1)}$ with a maximal compact subgroup $K_S = Pin(n-1) \times Pin(n-1)$.  Consider the embedding $\widetilde{i}:S \hookrightarrow \widetilde{Spin(n,n)_0} \hookrightarrow G$.   $\widetilde{i}(^0\!M_{\widetilde{Pin(n-1,n-1)}}) = {^0\!M_{G}}$ (cf. proof of Theorem 1 and Theorem 2.1).

Let $V_{\xi}$ be an irreducible $K$ module that occurs in $\mathcal{H} \otimes V_{\tau}$ and let $V_{\gamma_1},...,V_{\gamma_N}$ be the distinct $K$ types that occur in $\mathcal{H}$ such that $V_{\xi} \subset V_{\gamma_j} \otimes V_{\tau}$ for all $j$.  The assumption can be replaced with $V_{\gamma_j} \subset V_{\xi} \otimes V_{\tau}^{*}$.  Let Span $K_S \cdot V_{\gamma_j}^{^0\!M_{G}} = \oplus_{k} V_{\gamma_j,k}$ be a decomposition into irreducible $K_S$ modules where each $V_{\gamma_{j,k}}$ occurs in $\mathcal{H}$ for $S$, and let $V_{\xi}|_{K_S} = \oplus_{l} V_{\xi,l}$ be a decomposition into irreducible $K_S$ submodules.  $V_{\xi,l}$ and $V_{\gamma_{j,k}}$ are outer tensor products of irreducible $Pin(n-1)$ modules.  If $V_{\xi,l}$ and $V_{\gamma_{j,k}}$ are such that $V_{\gamma_{j,k}} \subset V_{\xi,l} \otimes V_{\tau}^{*}$, either the first or the second factors of the two must be isomorphic as $Pin(n-1)$ modules depending on $V_{\tau}$.  As $t_{\alpha}$ commutes with $Spin(n-1) \times Spin(n-1)$ and ${^0\!M_{G}}$ acts on $t_{\alpha}$ by $\pm 1$, there is a unique dominant $t_{\alpha}$ weight on an irreducible $K_S$ module.  $V_{\gamma_{j,k}} \subset V_{\xi,l} \otimes V_{\tau}^{*}$ implies $V_{\xi,l} \subset V_{\gamma_{j,k}}\otimes V_{\tau}$ hence $V_{\xi,l}|_{^0\!M_G} = \oplus_{j} V_{\tau_j}$ where $V_{\tau_j} \cong V_{\tau}$ as $^0\!M_G$ modules by Lemma 4.12.  This is important as the $t_{\alpha}$ weights on $V_{\xi}$ that we consider are the ones from $^0\!M_G$ submodules of $V_{\xi}$ isomorphic to $V_{\tau}|_{^0\!M_{G}}$.  As no two isomorphic copies of $V_{\gamma_{j,k}}$ can occur in $V_{\xi,l} \otimes V_{\tau}^{*}$ ($n=3$ is well understood.  For $n\geq 4$, cf. proof of Lemma 4.11), the comparison of $t_{\alpha}$ weights is proved by restriction to $K_S$ as there is a unique dominant $t_{\alpha}$ weight on an irreducible $K_S$ module.

Let now $G$ be of type $F_4$ with a maximal compact subgroup $K = Sp(3) \times SU(2)$.  Let $S = \widetilde{Spin(5,4)_0}$ with a maximal compact subgroup $K_S = Spin(5) \times Spin(4)$.  Let $\widetilde{i}:S \hookrightarrow G$ be the embedding with $\widetilde{i}(^0\!M_{S}) = {^0\!M_{G}}$ (cf. proof of Theorem 1).  The small $K$ type $\mathbb{C}^2 \circ p_2$ restricted to $K_S$ is one of the two small $K_S$ types $s \circ p_2$ considered above.  Let $V_{\xi}$ be an irreducible $K$ module that occurs in $\mathcal{H} \otimes V_{\tau}$ and let $V_{\gamma_1},...,V_{\gamma_N}$ be the distinct $K$ types that occur in $\mathcal{H}$ such that $V_{\xi} \subset V_{\gamma_j} \otimes V_{\tau}$ for all $j$.  The assumption can be replaced with $V_{\gamma_j} \subset V_{\xi} \otimes V_{\tau}^{*}$.  Let Span $K_S \cdot V_{\gamma_j}^{^0\!M_{G}} = \oplus_{k} V_{\gamma_j,k}$ be a decomposition into irreducible $K_S$ modules where each $V_{\gamma_{j,k}}$ occurs in $\mathcal{H}$ for $S$, and let $V_{\xi}|_{K_S} = \oplus_{l} V_{\xi,l}$ be a decomposition into irreducible $K_S$ submodules.  Note that no two isomorphic copies of $V_{\gamma_{j,k}}$ can occur in $V_{\xi,l} \otimes V_{\tau}^{*}$ (cf. proof of Lemma 4.11).  This is important as the statement of the Theorem for $B_4$ assumes distinct $V_{\gamma_{j,k}}$ types.  Hence the comparison of $t_{\alpha}$ weights is that for the type of $B_4$ on the modules $\{V_{\xi_l}\}$ and $\{ V_{\gamma_{j,k}} \}$ similarly as before.
\end{proof}
\end{theorem}

\end{subsection}

\begin{subsection}{}
Recall the notation of 4.1.  For $\alpha$ a positive root of $\mathfrak{g}_{\mathbb{R}}$, let $\mathfrak{g_{\alpha}} = \mathbb{C}\theta(e_{\alpha}) \oplus \mathfrak{a} \oplus \mathbb{C} e_{\alpha} = \theta(\mathfrak{n}_{\alpha}) \oplus \mathfrak{a} \oplus \mathfrak{n}_{\alpha}$ be the triangular decomposition and let $\mathfrak{g_{\alpha}} = \mathbb{C} (e_{\alpha} + \theta(e_{\alpha})) \oplus (\mathbb{C} (e_{\alpha} - \theta(e_{\alpha})) \oplus \mathfrak{a}) = \mathfrak{k}_{\alpha} \oplus \mathfrak{p}_{\alpha}$ be the Cartan decomposition.  Let $\mathcal{H}_{\alpha}$ be the space of harmonics on $\mathfrak{p}_{\alpha}$ and let $\mathcal{J}_{\alpha} = S(\mathfrak{p}_{\alpha})^{\mathfrak{k}_\alpha}$.  For $\alpha$ simple, let $\mathfrak{n}^{\alpha} = \bigoplus_{\psi \in \Phi^{+} - \{\alpha\}} \mathfrak{g}^{\psi}$.  Let $\mathfrak{k}^{\alpha} = \oplus_{\psi \in \Phi^{+} - \{\alpha\}} \mathbb{C}(e_{\psi}+\theta e_{\psi})$ so that $\mathfrak{k} = \mathfrak{k}_{\alpha} \oplus \mathfrak{k}^{\alpha}$.  $\mathfrak{g} = \mathfrak{n}^{\alpha} \oplus \mathfrak{g}_{\alpha} \oplus \mathfrak{k}^{\alpha}$ and $\mathfrak{n}^{\alpha}$ is a Lie subalgebra of $\mathfrak{g}$ as $\alpha$ is simple.

\begin{lemma} For $\alpha$ a simple root of $\mathfrak{g}_{\mathbb{R}}$, 
\[ U(\mathfrak{g}) = symm(\mathcal{H}_{\alpha}) symm(\mathcal{J}_{\alpha}) U(\mathfrak{k}) \oplus \mathfrak{n}^{\alpha} U(\mathfrak{n}^{\alpha})symm(\mathcal{H}_{\alpha}) symm(\mathcal{J}_{\alpha}) U(\mathfrak{k}) \]
\begin{proof}
By Proposition 2.4.1 of [Kos], 
\[ U(\mathfrak{g}) = U(\mathfrak{n}^{\alpha})symm(\mathcal{H}_{\alpha})symm(\mathcal{J}_{\alpha}) \oplus U(\mathfrak{g})\mathfrak{k}\]  
Hence, 
\begin{align*}
U(\mathfrak{g}) &= 
U(\mathfrak{n}^{\alpha})symm(\mathcal{H}_{\alpha})symm(\mathcal{J}_{\alpha})U(\mathfrak{k})\\ &=
symm(\mathcal{H}_{\alpha})symm(\mathcal{J}_{\alpha})U(\mathfrak{k}) \oplus \mathfrak{n}^{\alpha}U(\mathfrak{n}^{\alpha})symm(\mathcal{H}_{\alpha})symm(\mathcal{J}_{\alpha})U(\mathfrak{k})
\end{align*}
\end{proof}
\end{lemma}

As $U(\mathfrak{g}) \otimes_{U(\mathfrak{k})U(\mathfrak{g})^{\mathfrak{k}}} V_{\tau} \cong I_{P,\sigma,\nu}$ as $K$ modules (cf. 11.3.6 [RRG II]), we have the following $K$ module isomorphisms by Lemma 5.7.
\begin{align*}
I_{P,\sigma,\nu} &\cong U(\mathfrak{g}) \otimes_{U(\mathfrak{k})U(\mathfrak{g})^{\mathfrak{k}}} V_{\tau} \\\
&\cong symm(\mathcal{H}_{\alpha}) symm(\mathcal{J}_{\alpha}) \otimes_{U(\mathfrak{k})U(\mathfrak{g})^{\mathfrak{k}}} V_{\tau} \\\ 
&\ \ \ \ \oplus \ \mathfrak{n}^{\alpha} U(\mathfrak{n}^{\alpha})symm(\mathcal{H}_{\alpha})symm(\mathcal{J}_{\alpha}) \otimes_{U(\mathfrak{k})U(\mathfrak{g})^{\mathfrak{k}}} V_{\tau}
\end{align*}
Let $\mathbb{C}[{^0\!M}]$ be the group algebra generated by ${^0\!M}$.  Denote by $U(\mathfrak{k}_{\alpha}) \Pi \mathbb{C}[{^0\!M}]$ the smash product of $U(\mathfrak{k}_{\alpha})$ with $\mathbb{C}[{^0\!M}]$, i.e. $U(\mathfrak{k}_{\alpha}) \Pi \mathbb{C}[{^0\!M}]$ has a ($U(\mathfrak{k}_{\alpha}),\mathbb{C}[{^0\!M}]$) action on $symm(\mathcal{H}_{\alpha}) \otimes V_{\tau}$ that is an analog of a ($\mathfrak{g},K$) action.  Let $I_{\tau} = U(\mathfrak{k}) \cap ker \ \tau$.  As ${^0\!M}$ acts irreducibly on $V_{\tau}$, $U(\mathfrak{k})/I_{\tau} \cong End(V_{\tau}) \cong (U(\mathfrak{k}_{\alpha})\Pi \mathbb{C}[{^0\!M}])/(ker \ \tau \cap (U(\mathfrak{k}_{\alpha}) \Pi \mathbb{C}[{^0\!M}]))$.

For $\alpha$ simple, let \[L_{\alpha}:U(\mathfrak{g}) \otimes_{U(\mathfrak{k})U(\mathfrak{g})^{\mathfrak{k}}} V_{\tau} \rightarrow symm(\mathcal{H}_{\alpha}) symm(\mathcal{J}_{\alpha}) \otimes_{U(\mathfrak{k})U(\mathfrak{g})^{\mathfrak{k}}} V_{\tau}\]
\[ \ \ \ \ \ \ \ \ \ \ \ \ \ \ \ \ \ \ \ \ \ \ \ \ \ \ \ \ \ \ \ \oplus \ \mathfrak{n}^{\alpha} U(\mathfrak{n}^{\alpha})symm (\mathcal{H}_{\alpha})symm(\mathcal{J}_{\alpha}) \otimes_{U(\mathfrak{k})U(\mathfrak{g})^{\mathfrak{k}}} V_{\tau}\] be the projection onto the first summand.  Denote by $Q$ the projection onto the first summand of $U(\mathfrak{g}) = U(\mathfrak{a})U(\mathfrak{k}) \oplus \mathfrak{n}U(\mathfrak{g})$ followed by the projection onto $U(\mathfrak{a}) \otimes (U(\mathfrak{k})/I_{\tau}) = U(\mathfrak{a}) \otimes End(V_{\tau})$.

Let $\rho = \frac{1}{2} \Sigma_{\phi \in \Phi^{+}} \phi$.  Let $\alpha$ be a positive root of $\mathfrak{g}_{\mathbb{R}}$.  Let $V_{\xi}$ be an irreducible $K$ module that occurs in $\mathcal{H} \otimes V_{\tau}$ and let $V_{\xi} = \bigoplus_{j=1}^{n(\xi)} V_{\tau_{j}^{\alpha}} \oplus U_{\alpha}$ be a decomposition into $K_{\alpha}$ submodules where each $V_{\tau_{j}^{\alpha}}$ is an irreducible $K_{\alpha}$ module such that $V_{\tau_j}^{\alpha} \cong V_{\tau}$ as $^0\!M$ modules for $j=1,...,n(\xi)$.  Let $\epsilon_{1},...,\epsilon_{n(\xi)}$ be a basis of $Hom_{K}(V_{\xi},symm(\mathcal{H})\otimes V_{\tau}) \cong Hom_{K}(V_{\xi},I_{P,\sigma,\nu})$.  Let $p_{\tau_{j}^{\alpha}}$ be the determinant of $P^{\tau_{j}^{\alpha}}$ matrix for the rank one subgroup $G_{\alpha}$ with $K_{\alpha}$ type ${\tau_{j}^{\alpha}}$.

\begin{theorem} If $\alpha$ is a simple root of $\mathfrak{g}_{\mathbb{R}}$, $p_{\tau_{j}^{\alpha}} \ | \ p_{\xi}$.
\begin{proof}

This proof follows that of Proposition 2.4.3 of [Kos] closely.

$(P^{\xi})_{ij}$ is the action of $\epsilon_i(V_{\tau_j}^{\alpha})(e)$ followed by the replacement of the element in $\mathbb{C}[\nu]$ with the corresponding element in $S(\mathfrak{a})$ where $e \in G$ is the identity element.  Recall that $V_{\tau_{j}^{\alpha}}$ is an irreducible $K_{\alpha}$ module isomorphic to a submodule of $\mathcal{H}_{\alpha} \otimes V_{\tau}$ as $\mathcal{H}_{\alpha} \otimes V_{\tau} \cong Ind_{^0\!M}^{K_{\alpha}}(V_{\tau})$ as $K_{\alpha}$ modules.  $L_{\alpha}$ is a $K_{\alpha}$ map as $[\mathfrak{g}_{\alpha},\mathfrak{n}^{\alpha}]\subseteq \mathfrak{n}^{\alpha}$ for $\alpha$ simple.  Therefore, $(P^{\xi})_{ij}$ is achieved via $L_{\alpha}(\epsilon_{i}(V_{\tau_{j}^{\alpha}}))(e)$ as the subspace $\mathfrak{n}^{\alpha}$ acts trivially at the identity $e$.

First assume $t_{\alpha}$ acts nontrivially on $V_{\tau}$.  Recall the definition of $V_{\tau,\alpha}^{\pm}$ from 4.1.  Without loss of generality, assume by Theorem 4.7 $L_{\alpha}(\epsilon_{i}(V_{\tau_{j}^{\alpha}}^{+})) = \overline{Z_{\alpha}}^{l_{j}}R_{i,j}^{\alpha} \otimes V_{\tau,\alpha}^{+}$ where $V_{\tau_{j}^{\alpha}}^{+}$ is the subspace of $V_{\tau_j}^{\alpha}$ that consists of dominant $t_{\alpha}$ weight vectors with $R_{i,j}^{\alpha} \in symm(\mathcal{J}_{\alpha})$.  $symm(\mathcal{J}_{\alpha}) \subseteq U(\mathfrak{g}_{\alpha})^{\mathfrak{k}_{\alpha}}$ with $U(\mathfrak{g}_{\alpha})^{\mathfrak{k}_{\alpha}}$ the subalgebra generated by $t_{\alpha}$, the center of $\mathfrak{g}_{\alpha}$, and the Casimir element.  The action of $\overline{Z_{\alpha}}^{l_{j}}R_{i,j}^{\alpha}$ on $V_{\tau,\alpha}^{+}$ at $e$ is $Q(\overline{Z_{\alpha}}^{l_{j}}R_{i,j}^{\alpha})$.  By 3.5.6 of [RRG I], $Q(\overline{Z_{\alpha}}^{l_{j}}R_{i,j}^{\alpha}) = Q(R_{i,j}^{\alpha})Q(\overline{Z_{\alpha}}^{l_{j}}) = r_{i,j}^{\alpha}Q(\overline{Z_{\alpha}}^{l_{j}})$ where the action of $r_{i,j}^{\alpha}$ on $V_{\tau,\alpha}^{+}$ is invariant under $\tilde{x_{\alpha}} = T_{-\rho}x_{\alpha}T_{\rho}$ the translated Weyl group element of simple reflection.  $Q(\overline{Z_{\alpha}}^{l_{j}})$ is an element of $U(\mathfrak{a}) \otimes (U(\mathfrak{k})/I_{\tau})$ where $U(\mathfrak{a}) \otimes (U(\mathfrak{k})/I_{\tau}) = U(\mathfrak{a}) \otimes End(V_{\tau}) = U(\mathfrak{a}) \otimes (U(\mathfrak{k}_{\alpha})\Pi \mathbb{C}[{^0\!M}])/(ker \ \tau \cap (U(\mathfrak{k}_{\alpha}) \Pi \mathbb{C}[{^0\!M}]))$.  Hence the action of $Q(\overline{Z_{\alpha}}^{l_{j}})$ on $V_{\tau,\alpha}^{+}$ is the determinant $p_{\tau_{j}}^{\alpha}$ of $P^{\tau_{j}}$ matrix for the rank one subgroup $G_{\alpha}$ with $K_{\alpha}$ type $\tau_{j}$, and $p_{\tau_{j}}^{\alpha}$ divides $p_{\xi}$.

Now assume $t_{\alpha}$ acts trivially on $V_{\tau}$.  There is a unique dominant $t_{\alpha}$ weight on $V_{\tau_{j}^{\alpha}}$ by Corollary 4.8.  Let $v_j \in V_{\tau_{j}^{\alpha}}$ be a dominant $t_{\alpha}$ weight vector.  As irreducible $K_{\alpha}$ submodules of $\mathcal{H}_{\alpha} \otimes V_{\tau}$ occur with multiplicity one, $L_{\alpha}(\epsilon_i(v_j)) \in \overline{Z_{\alpha}^l} symm(\mathcal{J}_{\alpha}) \otimes V_{\tau}$ for some $l \in \mathbb{Z}_{\geq 0}$ where $symm(\mathcal{J}_{\alpha}) \subseteq U(\mathfrak{g}_{\alpha})^{\mathfrak{k}_{\alpha}}$ with $U(\mathfrak{g}_{\alpha})^{\mathfrak{k}_{\alpha}}$ the subalgebra generated by $t_{\alpha}$, center of $\mathfrak{g}_{\alpha}$, and the Casimir element.  As $t_{\alpha}$ acts trivially on $V_{\tau}$, the action of $\overline{Z_{\alpha}^l} symm(\mathcal{J}_{\alpha})$ on $V_{\tau}$ at $e$ is given by $Q(\overline{Z_{\alpha}^l} symm(\mathcal{J}_{\alpha}))$, hence the rest of the argument is exactly the above.
\end{proof}
\end{theorem}

Recall the notation of 4.1.  Let $\rho_{\phi}$ play the role of $\rho$ for the rank one subalgebra $\mathfrak{g}_{\phi}$ for a positive root $\phi$ of $\mathfrak{g}_{\mathbb{R}}$.  Define $p_{\phi} = p_{\tau_{1}^{\phi}} ... p_{\tau_{n(\xi)}^{\phi}}$ where $p_{\tau_{j}^{\phi}}$ denotes the determinant of $P^{\tau_{j}^{\phi}}$ matrix of the rank one subgroup $G_{\phi}$ with $K_{\phi}$ type ${\tau_{j}^{\phi}}$.  Define $p_{(\phi)} = T_{\rho_{\phi}-\rho}p_{\phi}$ where $T_{\rho_{\phi}-\rho}$ is the translation by $\rho_{\phi}-\rho$.  Each $p_{\tau_{j}^{\phi}}$ is a polynomial in $h_{\phi} \in \mathfrak{a}$.  Hence $p_{\phi}$ and $p_{(\phi)}$ are also.  As $T_{\rho_{\phi}-\rho}(h_{\phi})=h_{\phi}$ for $\phi$ simple, $p_{\phi} = p_{(\phi)}$ for $\phi$ simple.

\begin{theorem} $p_{(\phi)} \ | \ p_{\xi}$ for any $\phi \in \Phi^{+}$.
\begin{proof}
This proof is almost the same as that of Proposition 2.4.5 of [Kos].
			
Denote by $\Delta \subset \Phi^{+}$ the simple root system.  Let $\phi \in \Phi^{+}$ and let $O(\phi) = \Sigma m_i$ if $\phi = \Sigma_{\alpha_{i} \in \Delta}m_i\alpha_i$.  If $O(\phi)=1$, $\phi \in \Delta$ hence the claim is true by Theorem 5.8.  We proceed by induction on $O(\phi)$.  Assume $O(\phi)>1$ and assume the claim for all $\psi \in \Phi^{+}$ with $O(\psi)<O(\phi)$.  For some $\alpha \in \Delta$, $\left\langle{\phi,h_{\alpha}}\right\rangle > 0$ and let $\psi \in \Phi^{+}$ be such that $O(\psi)<O(\phi)$.  $\phi = x_{\alpha} \cdot \psi$ for some $\alpha \in \Delta$ where $x_{\alpha}$ is the Weyl group element of simple reflection.  Note $\psi \neq \alpha$.  $p_{\xi} = r_{\alpha} p_{\alpha}$ where $r_{\alpha}$ is invariant under the action of $\tilde{x_{\alpha}}$ (cf. proof of Theorem 5.8).  By induction hypothesis, $p_{(\psi)}$ divides $r_{\alpha} p_{\alpha}$.  As $p_{\alpha}$ is a polynomial in $h_{\alpha}$ whereas $p_{(\psi)}$ is a polynomial in $h_{\psi}$, and as $h_{\alpha}\neq h_{\psi}$, $p_{\alpha}$ and $p_{(\psi)}$ are mutually prime.  Hence $p_{(\psi)}$ divides $r_{\alpha}$, therefore $\tilde{x_{\alpha}} \cdot p_{(\psi)}$ does also.

We assert $\tilde{x_{\alpha}} \cdot p_{(\psi)} = p_{(\phi)}$ up to a nonzero scalar.  As $x_{\alpha} \cdot \psi=\phi$, $x_{\alpha} \cdot \mathfrak{g}_{\psi}=\mathfrak{g}_{\phi}$, $x_{\alpha}\cdot \mathfrak{k}_{\psi}=\mathfrak{k}_{\phi}$, and $x_{\alpha} \cdot \mathfrak{p}_{\psi}=\mathfrak{p}_{\phi}$.  Moreover, $x_{\alpha}\cdot \mathfrak{a}=\mathfrak{a}$ and $x_{\alpha}\cdot \mathfrak{n}_{\psi}=\mathfrak{n}_{\phi}$.  Therefore, for $u\in U(\mathfrak{g}_{\psi})$, $x_{\alpha} \cdot Q(u)=Q(x_{\alpha} \cdot u)$.  Also, $x_{\alpha}K_{\psi}x_{\alpha}^{-1}=K_{\phi}$.  Furthermore, let $V_{\xi}=\oplus_{j=1}^{n(\xi)} V_{\tau_{j}^{\psi}} \oplus U_{\psi}$ be a decomposition into $K_{\psi}$ submodules such that each $V_{\tau_{j}^{\psi}}$ is an irreducible $K_{\psi}$ module isomorphic to $V_{\tau}|_{^0\!M}$ when restricted to ${^0\!M}$.  If $V_{\tau_{j}^{\phi}}=x_{\alpha} \cdot V_{\tau_{j}^{\psi}}$ and $U_{\phi} = x_{\alpha} \cdot U_{\psi}$, then $V_{\xi}=\oplus_{j=1}^{n(\xi)} V_{\tau_{j}^{\phi}} \oplus U_{\phi}$ is a decomposition into $K_{\phi}$ submodules such that each $V_{\tau_{j}^{\phi}}$ is an irreducible $K_{\phi}$ module isomorphic to $V_{\tau}|_{^0\!M}$ when restricted to ${^0\!M}$ by Lemma 4.13.  Hence $x_{\alpha} \cdot p_{\psi}=p_{\phi}$ up to a nonzero scalar.  But, $\tilde{x_{\alpha}} \cdot p_{(\psi)}=T_{-\rho}x_{\alpha}T_{\rho}T_{\rho_{\psi}-\rho}p_{\psi}=T_{-\rho}x_{\alpha}T_{\rho_{\psi}}p_{\psi}=T_{-\rho}x_{\alpha}T_{\rho_{\psi}}x_{\alpha}^{-1}x_{\alpha} \cdot p_{\psi} = T_{-\rho}T_{x_{\alpha}\cdot \rho_{\psi}}x_{\alpha}\cdot p_{\psi}=T_{\rho_{\phi}-\rho}p_{\phi}=p_{(\phi)}$.  This completes the assertion and $p_{(\phi)}$ divides $p_{\xi}$.
		
\end{proof}
\end{theorem}
\end{subsection}

\begin{subsection}{}

Let $V_{\xi}$ be an irreducible $K$ module that occurs in $\mathcal{H} \otimes V_{\tau}$.  Recall the notation of section 3 and section 5.2.  We are now ready to prove Theorem 2 in the introduction.

\begin{proof} (of Theorem 2)

The right hand side divides the left hand side by Theorem 5.9 and the fact that for $\phi, \psi \in \Phi^{+}$ with $\phi \neq \psi$, $p_{(\phi)}$ and $p_{(\psi)}$ are mutually prime.  We assert the degrees of the two polynomials are equal.  As the definition of $p_{\xi}$ is independent of the basis up to a nonzero scalar, we may assume $\epsilon_{i}^{\xi}(T_{j}^{\xi}(V_{\tau})) \subseteq symm(\mathcal{H})_{deg(i)} \otimes V_{\tau}$, i.e. a homogeneous basis.  Therefore, the degree of the left hand side, $deg(\xi)$, is at most $\Sigma_{i=1}^{n(\xi)} deg(i)$.  If ${\gamma_{1}},...,{\gamma_{N}}$ are the distinct $K$ types that occur in $symm(\mathcal{H})$ with $V_{\xi} \subseteq  V_{\gamma_{j}}\otimes V_{\tau}$,  $\Sigma_{i=1}^{n(\xi)} deg(i) = \Sigma_{j=1}^{N} deg(\gamma_{j})$ where $deg(\gamma_{j})$ is the sum of the degrees in which $V_{\gamma_{j}}$ occur in $symm(\mathcal{H})$.  But, $deg(\gamma_{j})$ = $\Sigma_{\phi \in \Phi^{+}} n_{\gamma_{j}}^{\phi}$ where $n_{\gamma_{j}}^{\phi}$ is the sum of the degrees in which the irreducible $K_{\phi}$ submodules in Span $K_{\phi} \cdot V_{\gamma_{j}}^{{^0\!M}}$ occur in $symm(\mathcal{H}_{\phi})$ by Proposition 2.3.12 and Theorem 2.3.14 of [Kos].  $n_{\gamma_{j}}^{\phi}$ only depends on the unique dominant $t_{\phi}$ weights (cf. Theorem 4.6).  

$\Sigma_{i=1}^{n(\xi)} deg(i)$ = $\Sigma_{j=1}^{N}deg(\gamma_{j})$=$\Sigma_{j=1}^{N} \Sigma_{\phi \in \Phi^{+}} n_{\gamma_{j}}^{\phi}$= $\Sigma_{\phi \in \Phi^{+}} deg(p_{(\phi)}(\nu))$ by Theorem 5.4, Theorem 5.6, and the fact that $deg(Q(Z_{\phi}^{l}))=l$ by Theorem 7.6 of [JW].  Indeed, first 
consider all other than type $B_n$.  For $\phi$ such that $t_{\phi}$ acts nontrivially on $V_{\tau}$, an irreducible $K_{\phi}$ module that occurs in $\mathcal{H}_{\phi} \otimes V_{\tau}$ is uniquely determined by the unique, dominant $t_{\phi}$ weight (cf. Lemma 4.2, Lemma 4.3, Theorem 4.7).  For $\phi$ such that $t_{\phi}$ acts trivially on $V_{\tau}$, although the two irreducible $K_{\phi}$ modules that occur in $(\mathbb{C}Z_{\phi}^{l} \oplus \mathbb{C}\overline{Z_{\phi}^{l}}) \otimes V_{\tau}$ are not distinguishable by the unique, dominant $t_{\phi}$ weight, degree is the same for the two.  Now consider type $B_n$.  The argument for long roots is similar to the above as the long roots act by $\pm \frac{1}{2}$ on the small $K$ types (cf. Lemma 4.2).  For a short root $\phi$, $t_{\phi}$ acts trivially on the small $K$ type $s \circ p_2$ (cf. Lemma 4.3), hence the argument is again similar to the above.  For a short root $\phi$, $t_{\phi}$ acts by $\pm 1$ on the small $K$ type $s \circ p_1$.  In this case, an irreducible $K_{\phi}$ module that occurs in $\mathcal{H}_{\phi} \otimes V_{\tau}$ is not uniquely determined by the unique, dominant $t_{\phi}$ weight.  But, recall from the proof of Theorem 5.6 the irreducible $K_S = Pin(n-1) \times Pin(n-1)$ modules $V_{\xi,l}$, $V_{\gamma_{j,k}}$ such that $V_{\xi,l} \subset V_{\gamma_{j,k}}\otimes V_{\tau}$.  We may assume $\phi$ is such that $t_{\phi}$ has a unique dominant weight on $V_{\xi,l}$ and on $V_{\gamma_{j,k}}$.  If $v$ is a dominant $t_{\phi}$ weight vector of $V_{\xi,l}$, $L_{\phi}(v) \in symm(\mathcal{H}_{\phi})symm(\mathcal{J}_{\phi}) \otimes_{U(\mathfrak{k})U(\mathfrak{g})^{\mathfrak{k}}} V_{\tau}$ must have the factor in $symm(\mathcal{H}_{\phi})$ that corresponds to any of the irreducible $K_{\phi}$ submodules in $Span \ K_{\phi} \cdot V_{\gamma_{j,k}}^{^0\!M}$.  Indeed, $V_{\gamma_{j,k}}$ has a unique dominant $t_{\phi}$ weight, and $L_{\phi}$ is a $K_{\phi}$ map.

Therefore the degree of the left hand side is less than or equal to the degree of the right hand side.  As the right hand side divides the left hand side, the degree of the right hand side is less than or equal to the degree of the left hand side.  Hence the degrees of the two polynomials are equal.
\end{proof}

\end{subsection}

\end{section}
\begin{section}{Applications}

In this section, we derive implications and applications of the product formula for $p_{\xi}$.  Let $\mathfrak{g}_{\mathbb{R}}$ be any of the split real forms of simple Lie types other than type $C_n$.  Let $G$ be the connected, simply connected Lie group with Lie algebra $\mathfrak{g}_{\mathbb{R}}$.  In 6.1, the intertwining operator between the genuine principal series representations of $G$ is realized as a ratio of the $P^{\xi}$ matrices.  In 6.2, we compute $p_{\xi}$ explicitly for the rank one case of type $A_1$.  We use these results to derive a formula for $p_{\xi}$ and the determinant of the intertwining operator between the genuine principal series representations of $\widetilde{SL(n,\mathbb{R})}$.  In 6.3, we prove Theorem 3 stated in the introduction.  We note that in this section, the principal series representations, the intertwining operator between the principal series representations, and $P^{\xi}(\nu)$ are with $\rho$ shifts.  

\begin{subsection}{}

Let $(\rho_{\xi},V_{\xi})$ be an irreducible $K$ module that occurs in $\mathcal{H} \otimes V_{\tau}$ for a small $K$ type ${\tau}$.  Let $\sigma = V_{\tau}|_{^0\!M}$.  Recall the notation of section 3.  $P^{\xi}(\nu)$ is a map of $\bigoplus_{j=1}^{n(\xi)} T_{j}^{\xi}(V_{\tau}) \longrightarrow \bigoplus_{j=1}^{n(\xi)} T_{j}^{\xi}(V_{\tau})$ by setting $P^{\xi}(\nu)T_{i}^{\xi}(v) = \Sigma_{j=1}^{n(\xi)} T_{j}^{\xi}(P_{ji}^{\xi}(\nu)v)$ for $v \in V_{\tau}$.  Let $V_{\xi} = \bigoplus_{j=1}^{n(\xi)} T_{j}^{\xi}(V_{\tau}) \oplus U$ be a decomposition of $V_{\xi}$ as a ${}^{0}\!M$ module.  We define a map of $P^{\xi}(\nu)$ on $V_{\xi}$ where $P^{\xi}(\nu)$ acts as above on $\bigoplus_{j=1}^{n(\xi)} T_{j}^{\xi}(V_{\tau})$ and $P^{\xi}(\nu)$ acts by identity on $U$.  Thus, we may consider $P^{\xi}(\nu)$ as an operator on $Hom_{{}^{0}\!M}(V_{\xi},V_{\tau})$ where $P^{\xi}(\nu) \cdot \lambda = \lambda \circ P^{\xi}(\nu)$ for $\lambda \in Hom_{{}^{0}\!M}(V_{\xi},V_{\tau})$.

For $\lambda \in Hom_{^0\!M}(V_{\xi},V_{\tau})$ and $v \in V_{\xi}$, define ($\lambda \otimes v$)($k$) = $\lambda(\rho_{\xi}(k)v)$ that gives a $K$ module isomorphism $Hom_{^0\!M}(V_{\xi},V_{\tau}) \otimes V_{\xi} \cong I_{P,\sigma,\nu}(\xi)$ where $I_{P,\sigma,\nu}(\xi)$ is the $\xi$ isotypic component of $I_{P,\sigma,\nu}$.  Define for $a \in Hom_{K}(V_{\xi}, symm(\mathcal{H}) \otimes V_{\tau})$ and $v \in V_{\xi}$, $B^{\xi}_{\nu}(a)(v) = \pi_{P,\sigma,\nu}(a(v)) (e)$.  $B^{\xi}_{\nu}: Hom_{K}(V_{\xi}, symm(\mathcal{H}) \otimes V_{\tau}) \longrightarrow Hom_{^0\!M}(V_{\xi},V_{\tau})$.  Let $T_{\nu}: symm(\mathcal{H}) \otimes V_{\tau} \longrightarrow I_{P,\sigma,\nu}$ be defined by $T_{\nu}(\Sigma (u_{j} \otimes v_{j})) = \Sigma (\pi_{P,\sigma,\nu}(u_{j}) v_{j})$.  Let $\nu_{0} \in \mathfrak{a}^{*}$ be such that $T_{\nu_{0}}$ is a bijection (11.3.6 [RRG II]).

\begin{lemma} $T_{\nu} \circ T_{\nu_{0}}^{-1} (\lambda \otimes v)$ = $\lambda \circ P^{\xi}(\nu) \otimes v$ for $\lambda \in Hom_{^0\!M}(V_{\xi},V_{\tau})$ and $v \in V_{\xi}$.
\begin{proof}
This proof is almost the same as that of Lemma 7.3 of [JW].
		
Let $\hat{\delta}_{\xi} : Hom_{K}(V_{\xi}, symm(\mathcal{H}) \otimes V_{\tau}) \longrightarrow Hom_{^0\!M}(V_{\xi},V_{\tau})$ be defined so that $T_{\nu_{0}} (a(v)) = \hat{\delta}_{\xi} (a) \otimes v$ for $a \in Hom_{K}(V_{\xi}, symm(\mathcal{H}) \otimes V_{\tau})$ and $v \in V_{\xi}$.  By the above discussion, $B_{\nu_0}^{\xi}(a) = \hat{\delta}_{\xi}(a)$.  $T_{\nu} \circ T_{\nu_{0}}^{-1} (B_{\nu_{0}}^{\xi}(a) \otimes v) = B_{\nu}^{\xi}(a) \otimes v$.  But $B_{\nu}^{\xi}(a_{i})(T_{j}^{\xi}(V_{\tau})) = P^{\xi}_{ij}(\nu)$ where the basis $\left \{ a_{i}\right \}$ of $Hom_{^0\!M}(V_{\xi},V_{\tau})$ and $\{ T_{j}^{\xi}(V_{\tau})\}$ can be chosen so that $P^{\xi}(\nu_{0})$ is an identity matrix as $P^{\xi}(\nu_{0})$ is invertible.  Thus $B_{\nu}^{\xi}(a)=B_{\nu_{0}}^{\xi}(a) \circ P^{\xi}(\nu)$.  Therefore, if $B_{\nu_{0}}^{\xi}(a)=\lambda$, $T_{\nu} \circ T_{\nu_{0}}^{-1} (\lambda \otimes v)$ = $\lambda \circ P^{\xi}(\nu) \otimes v$.
\end{proof}
\end{lemma}

Let $\overline{P}$ be the opposite parabolic subgroup of a minimal parabolic subgroup $P$.  Define $\overline{T}_{\nu}: symm(\mathcal{H}) \otimes V_{\tau} \longrightarrow I_{\overline{P},\sigma,\nu}$ as we did for $P$.

\begin{theorem} Let $A(\nu): I_{P,\sigma,\nu} \longrightarrow I_{\overline{P},\sigma ,\nu}$ be such that $A(\nu) w = w$ for $w \in V_{\tau}$ and $A(\nu ) \circ \pi_{P,\sigma,\nu}(u) = \pi_{\overline{P},\sigma,\nu} (u) \circ A(\nu)$ for all $u \in U(\mathfrak g)$.  Then 
\[ A(\nu) (\lambda \otimes v) = \lambda \circ P^{\xi}(\nu)^{-1} P^{\xi}(-\nu) \otimes v \]
for $\lambda \in Hom_{^0\!M}(V_{\xi},V_{\tau})$, $v \in V_{\xi}$, and $\nu$ such that $det \ P^{\xi}(\nu) \neq 0$.
\begin{proof}
This proof is almost the same as that of Lemma 7.5 of [JW].
	
If $u \otimes w \in symm(\mathcal{H}) \otimes V_{\tau}$ is a simple tensor, 
\begin{align*}
A(\nu)T_{\nu}(u\otimes w) &= 
A(\nu) \pi_{P,\sigma,\nu}(u)w \\ &= 
\pi_{\overline{P},\sigma, \nu}(u) A(\nu) w\\ &= 
\pi_{\overline{P},\sigma, \nu}(u) w \\ &= 
\overline{T}_{\nu}(u \otimes w)
\end{align*}
Hence $A(\nu) (T_{\nu} \circ T_{\nu_{0}}^{-1})(\lambda \otimes v) = (\overline{T}_{\nu} \circ \overline{T}_{-\nu_0}^{-1})(\lambda \otimes v)$.  Thus,
\[ A(\nu )(\lambda \circ P^{\xi}(\nu) \otimes v) = \lambda \circ P^{\xi}(-\nu) \otimes v \]
by Lemma 6.1.  Therefore,	

\[ A(\nu) (\lambda \otimes v) = \lambda \circ P^{\xi}(\nu )^{-1} P^{\xi}(-\nu) \otimes v \]
\end{proof}
\end{theorem}

\end{subsection}

\begin{subsection}{}

Let $G$ be any of the connected, simply connected split real forms of simple Lie types other than type $C_n$.  Recall the notation of 4.1.  For a positive root $\alpha$ of $\mathfrak{g}_{\mathbb{R}}$, ${\mathfrak{g}_{\mathbb{R}}}_{\alpha} \cong \mathfrak{sl}(2,\mathbb{R}) \oplus Z({\mathfrak{g}_{\mathbb{R}}}_{\alpha})$.  Recall for $G_{\alpha}$, \[ \mathcal{H}_{\alpha} = \bigoplus_{l\geq 0} Z_{\alpha}^l \oplus \bigoplus_{l>0} \overline{Z_{\alpha}^l} \] where $Z_{\alpha} = h_{\alpha} + iy_{\alpha}$.  We drop the subscript $\alpha$ for convenience.  Recall the projection map $Q:U(\mathfrak{g}) \longrightarrow U(\mathfrak{a})U(\mathfrak{k}) \oplus \mathfrak{n}U(\mathfrak{g})$ onto the first summand. To compute $p_{\xi}(\nu)$ for a $K$ type ${\xi}$ that occurs in $\mathcal{H} \otimes V_{\tau}$, we compute $Q(Z^l)$ and $Q(\overline{Z}^l)$ with the $t$ weights on $V_{\tau}$ (cf. Lemma 4.2 and Lemma 4.3).  To do this, we use $Q'(Z^l)$ and $Q'(\overline{Z}^l)$ already computed in [JW] where $Q':U(\mathfrak{g}) \longrightarrow U(\mathfrak{a}) \oplus \mathfrak{n}U(\mathfrak{g}) \oplus U(\mathfrak{g})\mathfrak{k}$ is the projection onto the first summand.  This is because $Q(Z^l)$ and $Q(\overline{Z}^l)$ can be written as a sum of two different parts, one in $U(\mathfrak{a})$, and the other in $U(\mathfrak{a})U(\mathfrak{k})\mathfrak{k}$, where the former is exactly $Q'(Z^l)$ and $Q'(\overline{Z}^l)$ respectively.

\begin{theorem} $Q(Z^l)$ = $\Pi_{j=0}^{l-1} (h+2j-t)$ and $Q(\overline{Z}^l)$ = $\Pi_{j=0}^{l-1} (h+2j+t)$.
\begin{proof}
We prove the first formula.

By Theorem 7.6 of [JW], $Q'(Z^l)$ = $\Pi_{j=0}^{l-1} (h+2j)$.  We find the shift from $U(\mathfrak{a})U(\mathfrak{k})\mathfrak{k}$ part as it is the only difference between $Q$ and $Q'$.  If $l=1$, then $Z= h+iy$ = $h+i(2e+i{t})$, hence the formula is true.  We proceed by induction.  Assume the formula for $l-1$.  We prove the formula for $l$.  $Z^{l} =  ZZ^{l-1} = (h+2ie-{t})Z^{l-1}$.  After dropping the $\mathfrak{n}$ part $e$, $(h-{t})Z^{l-1}$ is left.  There are exactly two $U(\mathfrak{a})U(\mathfrak{k})\mathfrak{k}$ shifts, one from $h(U(\mathfrak{a})U(\mathfrak{k})\mathfrak{k}$ part of $Z^{l-1}$) and the other from $-{t}(Z^{l-1}) = -Z^{l-1}({t}-2(l-1))$ by the commutation relation.  Hence, the overall shift is $(h+2(l-1))(Q(Z^{l-1})-Q'(Z^{l-1}))-tQ(Z^{l-1})$.  As $Q'(Z^l) = (h+2(l-1))Q'(Z^{l-1})$,
\begin{align*}
Q(Z^l)&= Q'(Z^l) + shift \\ &= 
(h+2(l-1))Q'(Z^{l-1}) \\ &
\ \ \ \ \ +(h+2(l-1))(Q(Z^{l-1})-Q'(Z^{l-1}))-tQ(Z^{l-1}) \\ &= 
(h+2(l-1))Q(Z^{l-1})-tQ(Z^{l-1})\\ &= 
(h+2(l-1)-t)Q(Z^{l-1})
\end{align*}
Hence the first formula is proved.  The second formula is proved similarly.
\end{proof}
\end{theorem}

Assume $t$ acts nontrivially on $V_{\tau}$ with dominant weight $r$.

For the $\xi$ type $\overline{Z}^{l} \otimes V_{\tau}^{+} \oplus Z^{l}\otimes V_{\tau}^{-}$, 
				\[ p_{\xi}(\nu) = \Pi_{j=0}^{l-1} (\nu + 2j + 1 + r) \] and for the $\xi$ type $\overline{Z}^{l} \otimes V_{\tau}^{-} \oplus Z^{l}\otimes V_{\tau}^{+}$, 
				\[ p_{\xi}(\nu) = \Pi_{j=0}^{l-1} (\nu + 2j + 1 - r) \]
				
Assume $t$ acts trivially on $V_{\tau}$.

For any that occurs in $(\mathbb{C}\overline{Z}^{l} \oplus \mathbb{C}Z^{l}) \otimes V_{\tau}$, 
				\[ p_{\xi}(\nu) = \Pi_{j=0}^{l-1} (\nu + 2j + 1) \]

Consider the group $\widetilde{SL(n,\mathbb{R})}$.  Define $q_{\nu}:2\mathbb{N}+1 \longrightarrow \mathbb{C}[\nu]$ as follows.
	
		$q_{\nu}(m) := \Pi_{l=0}^{\frac{m-1}{4}} \Pi_{j=0}^{l-1} (\nu+2j+\frac{1}{2}) (\nu+2j+\frac{3}{2})$ if $4 \ | \ m-1$
		
		$q_{\nu}(m) := \Pi_{l=0}^{\frac{m-3}{4}} \Pi_{j=0}^{l-1} (\nu+2j+\frac{1}{2}) (\nu+2j+\frac{3}{2})$ 

			$ \ {} \ {} \ {} \ {} \ {} \ {} \ {} \  \ {} \ {} \ {} \ {} \ {} \ {} \ {} \ \times \ \Pi_{k=0}^{\frac{m-3}{4}} (\nu+2k+\frac{1}{2})$ if $4 \ | \ m-3$
			
\vspace{.1 in}
		
Define $G_{\nu}(m):2\mathbb{N}+1 \longrightarrow M$ where $M$ is the space of meromorphic functions, as follows.
		
		$\Gamma_{\nu}(m):= \Pi_{l=0}^{\frac{m-1}{4}} \Pi_{j=0}^{l-1}  \frac{\Gamma (\nu-2j+\frac{1}{2})}{\Gamma (\nu-2j-\frac{3}{2})} \frac{\Gamma (\nu+2j+\frac{1}{2})}{\Gamma (\nu+2j+\frac{5}{2})} $ if $4 \ | \ m-1$
		
		$\Gamma_{\nu}(m):= \Pi_{l=0}^{\frac{m-3}{4}} \Pi_{j=0}^{l-1} \frac{\Gamma (\nu-2j+\frac{1}{2})}{\Gamma (\nu-2j-\frac{3}{2})} \frac{\Gamma (\nu+2j+\frac{1}{2})}{\Gamma (\nu+2j+\frac{5}{2})} \times \ \Pi_{k=0}^{\frac{m-3}{4}} \frac{\Gamma (\nu-2k+\frac{1}{2})}{\Gamma (\nu-2k-\frac{3}{2})} $ if $4 \ | \ m-3$
		
\vspace{.1 in}
		
Let $V_{\xi}$ be an irreducible $K=Spin(n)$ module that occurs in $I_{P,\sigma,\nu}$.  The dominant $t_{\alpha}$ weights of interest on $V_{\xi}$ counting multiplicity are independent of $\alpha \in \Phi^{+}$.  Hence, we find the dominant $t_{\alpha}$ weights of interest on $V_{\xi}$ for $\alpha =  \epsilon_1-\epsilon_2$.  We branch down $V_{\xi}$ to $Spin(3)$ instead of $Spin(2) \cong SO(2)$ to simplify notation, where $Spin(3)$ is such that the quotient group $SO(3)$ occurs in the top left corner of the quotient group $SO(n)$ of $Spin(n)$.  Let $n = 2k+1$ or $n=2k$ and let $\xi = \xi_1\epsilon_1 + ... + \xi_k\epsilon_k$ be the highest weight of $\xi$.  Let $\frac{j_1}{2},...,\frac{j_{m_{\xi}}}{2}$ be the highest weights of irreducible $Spin(3)$ modules that occur in the branching counting multiplicity.

\[ p_{\xi}(\nu) = (\Pi_{\alpha \in \Phi^{+}} \Pi_{k=1}^{m_{\xi}} q_{(\nu,\alpha)}(j_k))^{\frac{2}{dim(V_{\tau})}} \]
\begin{align*}
det A(\nu)|_{I_{P,\tau,\nu}(\xi)} &= 
(\frac{p_{\xi}(-\nu)}{p_{\xi}(\nu)})^{dim(V_{\xi})}\\ &=  
((\Pi_{\alpha \in \Phi^{+}} \Pi_{k=1}^{m_{\xi}} \Gamma_{(\nu,\alpha)}(j_k))^{\frac{2}{dim(V_{\tau})}})^{dim(V_{\xi})} 
\end{align*}
	
$\noindent$where if $n=2k+1$, $dim(V_{\xi}) = \Pi_{1\leq i<j\leq k} \frac{(\xi_i+\rho_i)^2-(\xi_j+\rho_j)^2}{\rho_i^2-\rho_j^2} \Pi_{1\leq i\leq k} \frac{\xi_i+\rho_i}{\rho_i}$ with $\rho_i = k-i+\frac{1}{2}$ (cf. 7.1.2 [GW]), $dim(V_{\tau})=2^{k}$, and if $n=2k$, $dim(V_{\xi}) = \Pi_{1\leq i<j\leq k} \frac{(\xi_i+\rho_i)^2-(\xi_j+\rho_j)^2}{\rho_i^2-\rho_j^2}$ with $\rho_i=k-i$ (cf. 7.1.2 [GW]), $dim(V_{\tau})=2^{k-1}$.
\end{subsection}

\begin{subsection}{}

Consider the $(\mathfrak{g},K)$ module homomorphism
\[ \mu_{\tau,\nu}: U(\mathfrak{g}) \otimes_{U(\mathfrak{g})^{\mathfrak{k}}U(\mathfrak{k})} V_{\tau,\nu} \longrightarrow I_{P,\sigma,\nu} \] from 11.3.6 of [RRG II] where the homomorphism is the action of the first factor on the second as differential operators.  $U(\mathfrak{g}) \otimes_{U(\mathfrak{g})^{\mathfrak{k}}U(\mathfrak{k})} V_{\tau,\nu} \cong I_{P,\sigma,\nu}$ as $K$ modules (11.3.6 [RRG II]).  Therefore, the small $K$ type $I_{P,\sigma,\nu}(\tau)$ is cyclic if and only if $\mu_{\tau,\nu}$ is a $(\mathfrak{g},K)$ module isomorphism.  By the definition of $P^{\xi}(\nu)$, $I_{P,\sigma,\nu}(\tau)$ is cyclic if and only if $p_{\xi}(\nu) \neq 0$ for every $\xi \in \hat{K}$ that occurs in $I_{P,\sigma,\nu}$.

\begin{theorem} Let $G$ be any of the connected, simply connected, split real forms of simple Lie types other than type $C_n$.  

\begin{itemize}
\item Let $\tau$ be a small $K$ type other than $s\circ p_1$ for type $B_n$ ($n$ necessarily odd) and $\mathbb{C}^2\circ p_2$ for type $G_2$.  If $Re(\nu,\alpha) \geq 0$ for every $\alpha \in \Phi^{+}$, i.e. in the closed Langlands chamber, $I_{P,\sigma,\nu}(\tau)$ is cyclic.
\item Let $\tau = s\circ p_1$ for type $B_n$ ($n$ necessarily odd).  In the closed Langlands chamber, $I_{P,\sigma,\nu}(\tau)$ is cyclic for all other than $\nu$ such that $\frac{2(\nu,\alpha)}{(\alpha,\alpha)}=0$ for some $\alpha \in \Phi^{+}$ short, in which case is not cyclic.
\item Let $\tau = \mathbb{C}^2\circ p_2$ for type $G_2$.  In the closed Langlands chamber, $I_{P,\sigma,\nu}(\tau)$ is cyclic for all other than $\nu$ such that $\frac{2(\nu,\alpha)}{(\alpha,\alpha)}=\frac{1}{2}$ for some $\alpha \in \Phi^{+}$ short, in which case is not cyclic.
\end{itemize}
\begin{proof}
By Theorem 2, $p_{\xi}(\nu)$ up to a nonzero scalar is the product of those of rank one subgroups $G_{\alpha}$ of $G$ where $\alpha \in \Phi^{+}$.  As $G$ is split, the semisimple part of ${\mathfrak{g}_{\mathbb{R}}}_{\alpha}$ is isomorphic to $\mathfrak{sl}(2,\mathbb{R})$.  Let $n(\xi) = dim \ Hom_{^0\!M}(V_{\xi},V_{\tau})$.  

First consider the cases other than $s\circ p_1$ for type $B_n$ ($n$ necessarily odd) and $\mathbb{C}^2\circ p_2$ for type $G_2$.  The product formula for $p_{\xi}$, the formulas in the rank one case, Lemma 4.2, and Lemma 4.3 imply that for a $K$ type $\xi$ that occurs in $I_{P,\sigma,\nu}$, $\alpha \in \Phi^{+}$, $p_{\xi}(\nu)$ up to a nonzero scalar must be a product of polynomials of the form $(\frac{2(\nu,\alpha)}{(\alpha,\alpha )} + 2p + 1 \pm r)$ where $\alpha \in \Phi^{+}$, $r= 0, \ \frac{1}{2}$, and $p \in \mathbb{Z}_{\geq 0}$.  Therefore, if $Re(\nu,\alpha) \geq 0$ for every $\alpha \in \Phi^{+}$, $p_{\xi}(\nu) \neq 0$ for every $K$ type ${\xi}$ that occurs in $I_{P,\sigma,\nu}$ and $I_{P,\sigma,\nu}(\tau)$ is cyclic.

For type $B_n$, ${^0\!M} \cong {^0\!M_{\widetilde{SL(n,\mathbb{R})}}} \times \mathbb{Z}/2\mathbb{Z}$ where the $\mathbb{Z}/2\mathbb{Z}$ can be chosen to act trivially on $\mathcal{H}$ and on the small $K$ types $\tau = s \circ p_1$ ($n$ necessarily odd).  Therefore, the restriction of $\mathcal{H} \otimes V_{\tau}$ to $^0\!M$ just consists of $V_{\tau}|_{^0\!M}$ as there is only one genuine, irreducible representation of ${^0\!M_{\widetilde{SL(n,\mathbb{R})}}}$ ($n$ odd) up to isomorphism.  The product formula for $p_{\xi}$, the formulas in the rank one case, and Lemma 4.3 imply that there must be a $K$ type ${\xi}$ that occurs in $I_{P,\sigma,\nu}$ such that $p_{\xi}(\nu)$ contains factors of form $\frac{2(\nu,\alpha)}{(\alpha,\alpha )}$ for $\alpha \in \Phi^{+}$ short.  Therefore, $p_{\xi}(\nu)=0$ for $\nu$ such that $\frac{2(\nu,\alpha)}{(\alpha,\alpha)}=0$ for some $\alpha \in \Phi^{+}$ short.  For all other $\nu$ in the closed Langlands chamber, $p_{\xi}(\nu)\neq 0$ similarly as above.

As ${^0\!M} \cong {^0\!M_{\widetilde{SL(3,\mathbb{R})}}}$ for type $G_2$, there is only one genuine, irreducible representation of $^0\!M$ up
to isomorphism.  The product formula for $p_{\xi}$, the formulas in the rank one case, and Lemma 4.3 imply that there must be a $K$ type ${\xi}$ that occurs in $I_{P,\sigma,\nu}$ such that $p_{\xi}(\nu)$ contains factors of form $(\frac{2(\nu,\alpha)}{(\alpha,\alpha )} - \frac{1}{2})$ for $\alpha \in \Phi^{+}$ short.  Therefore, $p_{\xi}(\nu)=0$ for $\nu$ such that $\frac{2(\nu,\alpha)}{(\alpha,\alpha)}=\frac{1}{2}$ for some $\alpha \in \Phi^{+}$ short.  For all other $\nu$ in the closed Langlands chamber, $p_{\xi}(\nu)\neq 0$ similarly as above.
\end{proof}
\end{theorem}

\begin{theorem} Let $G$ be any of the connected, simply connected, split real forms of simple Lie types other than type $C_n$.  

\begin{itemize}
\item Let $\tau$ be a small $K$ type other than $s\circ p_1$ for type $B_n$ ($n$ necessarily odd).  Then, the underlying $(\mathfrak{g},K)$ module $I_{P,\sigma,\nu}$ $(Re \ \nu=0)$ of the unitary principal series is irreducible.
\item Let $\tau = s\circ p_1$ for type $B_n$ ($n$ necessarily odd).  Then, the underlying $(\mathfrak{g},K)$ module $I_{P,\sigma,\nu}$ $(Re \ \nu=0)$ of the unitary principal series is irreducible for all other than $\nu$ such that $\frac{2(\nu,\alpha)}{(\alpha,\alpha)}=0$ for some $\alpha \in \Phi^{+}$ short, in which case is reducible.
\end{itemize}
\begin{proof}
Consider first the case $\tau = s\circ p_1$ for type $B_n$ ($n$ necessarily odd) with $\nu$ is such that $\frac{2(\nu,\alpha)}{(\alpha,\alpha)}=0$ for some $\alpha \in \Phi^{+}$ short.  $I_{P,\sigma,\nu}(\tau)$ is not cyclic by Theorem 6.4, hence $I_{P,\sigma,\nu}$ is reducible.

Now consider all other cases.  If $I_{P,\sigma,\nu}$ is reducible, there is a proper, nontrivial, closed $(\mathfrak{g},K)$ invariant subspace $W$ of $I_{P,\sigma,\nu}$ that does not contain $\tau$ as $I_{P,\sigma,\nu}(\tau)$ is cyclic by Theorem 6.4.  Unitarity implies that the orthogonal complement of $W$, $W^{\perp}$, is a proper, nontrivial, closed $(\mathfrak{g},K)$ invariant subspace that contains $\tau$, which is a contradiction as $I_{P,\sigma,\nu}(\tau)$ is cyclic.
\end{proof}
\end{theorem}

We assert that there is a bijection between the set of equivalence classes of irreducible $(\mathfrak{g},K)$ modules admitting a small $K$ type ${\tau}$ and the set of Weyl group orbits in $\mathfrak{a}^{*}$ as in the spherical case.  Indeed, let $V$ be an irreducible $(\mathfrak{g},K)$ module admitting a small $K$ type ${\tau}$.  As $V_{\tau}|_{^0\!M}$ is irreducible, $V$ must be $(\mathfrak{g},K)$ isomorphic to a subquotient of $I_{P,\sigma,\nu}$ for some $\nu \in \mathfrak{a}^{*}$ by the Harish-Chandra Subquotient Theorem (cf. Theorem 3.5.6 [RRG I]).  Therefore, $V$ is completely determined by the action of $U(\mathfrak{g})^K$ on $V(\tau)$ by Proposition 3.5.4 of [RRG I].  Let $Q$ be the projection onto the first summand of $U(\mathfrak{g}) = U(\mathfrak{a})U(\mathfrak{k}) \oplus \mathfrak{n}U(\mathfrak{g})$ and let $\varpi(X) = X + \rho(X)$ for $X \in \mathfrak{a_{\mathbb{R}}}^{*}$.  The action of $U(\mathfrak{g})^K$ on $I_{P,\sigma,\nu}(\tau)$ is given by $((\nu \circ \varpi) \otimes \tau) (Q(g))$ for $g \in U(\mathfrak{g})^K$.  There exists a map $\gamma_{\tau}: U(\mathfrak{g})^K \rightarrow U(\mathfrak{a})^{N_{K}(A)/Z_K(A)}$ such that $(\varpi \otimes \tau) (Q(g)) = \gamma_{\tau}(g)$ (cf. Lemma 11.3.2 of [RRG II]).  In fact, $\gamma_{\tau}$ gives an algebra isomorphism $U(\mathfrak{g})^K /(U(\mathfrak{g})^K \cap U(\mathfrak{g})J_{\tau}) \cong U(\mathfrak{a})^{N_{K}(A)/Z_K(A)}$ where $J_{\tau} = ker \ \tau \subset U(\mathfrak{k})$.  If $\nu, \nu^{'} \in \mathfrak{a}^{*}$ and $\nu(f) = \nu^{'}(f)$ for all $f \in U(\mathfrak{a})^{N_{K}(A)/Z_K(A)}$, there exists $s\in N_{K}(A)/Z_K(A)$ such that $s\nu = \nu^{'}$ (cf. proof of Theorem 3.1.2 [RRG I]).  We are now ready to prove Theorem 3 in the introduction.

\begin{proof} (of Theorem 3)

$Y^{\tau,\nu}$ has a unique irreducible quotient (cf. 3.5.4 [RRG I]).  Thus $V$ is $(\mathfrak{g},K)$ isomorphic with the unique irreducible quotient of $Y^{\tau,\nu}$ for some $\nu$ in the closed Langlands chamber by the above discussion.  If $\tau$ is any other than $s \circ p_1$ for type $B_n$ ($n$ necessarily odd) and $\mathbb{C}^2 \circ p_2$ for type $G_2$, $\mu_{\tau,\nu}: U(\mathfrak{g}) \otimes_{U(\mathfrak{g})^{\mathfrak{k}}U(\mathfrak{k})} V_{\tau,\nu} \longrightarrow I_{P,\sigma,\nu}$ is a $(\mathfrak{g},K)$ isomorphism in the closed Langlands chamber by Theorem 6.4.
\end{proof}

Now consider all the cases other than $\tau = s \circ p_1$ for type $B_n$ ($n$ necessarily odd) and $\tau = \mathbb{C}^2 \circ p_2$ for type $G_2$.  Let $V$ be an irreducible $(\mathfrak{g},K)$ module that admits a small $K$ type ${\tau}$.  By Theorem 3, $V$ is $(\mathfrak{g},K)$ isomorphic with the unique irreducible quotient of $I_{P,\sigma,\nu}$ for some $\nu \in \mathfrak{a}^{*}$ in the closed Langlands chamber.  Suppose $Re(\nu,\alpha) > 0$ for some $\alpha \in \Phi^{+}$.  Let $F = \{\alpha \in \Delta \ | \ Re \ (\nu,\alpha) = 0\}$ and let ${\mathfrak{a}_{\mathbb{R}}}_{F} = \{ H \in \mathfrak{a}_{\mathbb{R}} \ | \ \alpha(H) = 0 \ \rm{for \ all} \ \alpha \in F\}$.  Let $^0\!M_{F} = Z_G({\mathfrak{a}_{\mathbb{R}}}_{F})$ and let $P_{F} = {^0\!M_{F}}A_F N_F$ be a given Langlands decomposition of a parabolic subgroup $P_F$ with special vector subgroup $A_F = exp ({\mathfrak{a}_{\mathbb{R}}}_{F})$.  Let $\nu = \varsigma + \mu $ where $\varsigma(H) = 0$, $\mu(H) = \nu(H)$ for $H \in {\mathfrak{a}_{\mathbb{R}}}_{F}$ and $\varsigma(H) = \nu(H)$, $\mu(H) = 0$ for $H \in \mathfrak{a}_{\mathbb{R}} - {\mathfrak{a}_{\mathbb{R}}}_{F}$.  $I_{P,\sigma,\nu}$ is realized in induction steps as $I_{P_F,{I_{{^0\!M_{F}} \cap P,\sigma,\varsigma}},\mu}$ where $(\sigma_F,{I_{{^0\!M_{F}} \cap P,\sigma,\varsigma}})$ is a tempered, unitary representation of $^0\!M_F$ and $\mu$ is in the Langlands chamber.  $(\sigma_F,{I_{{^0\!M_{F}} \cap P,\sigma,\varsigma}})$ is irreducible as $Ind$ is an exact functor and there exists an irreducible, unitary principal series representation $I_{P,\sigma,\chi}$ that is realized in induction steps as $I_{P_F,{I_{{^0\!M_{F}} \cap P,\sigma,\varsigma}},\chi^{'}}$ for some $\chi^{'} \in \mathfrak{a}^{*}$.  Otherwise, $Re(\nu,\alpha) = 0$ for all $\alpha \in \Phi^{+}$.  $I_{P,\sigma,\nu}$ then is a tempered, unitary principal series which is irreducible by Theorem 6.5.

\begin{corollary} Langlands parameters for $V$ are:
\begin{itemize}
\item $(P_F,\sigma_F,\mu)$ if $Re(\nu,\alpha) > 0$ for some $\alpha \in \Phi^{+}$.
\item Tempered if $Re(\nu,\alpha) = 0$ for all $\alpha \in \Phi^{+}$.
\end{itemize}
\end{corollary}

\begin{corollary} $V$ cannot be equivalent to a discrete series representation.
\begin{proof}
Let $F$ be as above.  If $F \neq \Delta$, then $V$ is not tempered.  If $F = \Delta$, then $V$ is equivalent with a unitary principal series representation for a proper parabolic subgroup which cannot be square integrable.
\end{proof}
\end{corollary}

\end{subsection}

\end{section}

Department of Mathematics, University of California, San Diego, La Jolla, CA, 92093.  E-mail: swl006@math.ucsd.edu

\end{document}